%
%  Quantum projective spaces and related q-hypergeometric
%  orthogonal polynomials
%
%  By M.S.Dijkhuizen and M. Noumi
%
%  email msdz@math.s.kobe-u.ac.jp
%        noumi@math.s.kobe-u.ac.jp
%
%  Version 16 April 1996
%  Last changes 13 May 1996
%
%  AMS-TeX  version 2.1
%
\input amstex
\font\thstil=cmsl10

\def\ltextindent#1{\hbox to x
\hangindent{#1\hss}\ignorespaces}

\newskip\proofskipamount
\proofskipamount=8pt plus 2pt minus 2 pt
\def\proofskip{\vskip\proofskipamount}
\def\endth{\par\ifdim\lastskip<\bigskipamount
\removelastskip\penalty55
\bigskip\fi}
\def\beginproof{\removelastskip
\penalty55\proofskip\noindent}
\def\endproof{\bigbreak}

\def\ifundefined#1{\expandafter\ifx
\csname#1\endcsname\relax}

\newif\ifdevelop \developfalse

\newtoks\chnumber
\newtoks\sectionnumber
\newcount\equationnumber
\newcount\thnumber

\def\assignnumber#1#2{%
\ifundefined{#1}\relax\else\message{#1 already defined}\fi
\expandafter\xdef\csname#1\endcsname
{\if-\the\sectionnumber\else\the\sectionnumber.\fi\the#2}%
}%

\def\beginsektion #1 #2 {\vskip0pt plus.1\vsize\penalty-250
\vskip0pt plus-.1\vsize
\bigbreak\bigskip
\sectionnumber{#1} \equationnumber0\thnumber0
\noindent{\bf #1. #2}\par
\nobreak\medskip\noindent}

\def\nologo{\expandafter\let
\csname logo\string @\endcsname=\empty}
\def\:{:\allowbreak }
\def\eq#1{\relax
\global\advance\equationnumber by 1
\assignnumber{EN#1}\equationnumber
{\tenrm (\csname EN#1\endcsname)}
}

\def\eqtag#1{\relax\ifundefined{EN#1}
\message{EN#1 undefined}{\sl (#1)}%
\else(\csname EN#1\endcsname)\fi%
}

\def\thname#1{\relax
\global\advance\thnumber by 1
\assignnumber{TH#1}\thnumber
\csname TH#1\endcsname
}

\def\beginth#1 #2 {\bigbreak\noindent{\bf #1\enspace \thname{#2}}
    ---\hskip4pt}
\def\thtag#1{\relax\ifundefined{TH#1}
\message{TH#1 undefined}{\sl #1}%
  \else\csname TH#1\endcsname\fi}

\def\al{\alpha}
\def\be{\beta}

\def\de{\delta}

\def\ze{\zeta}

\def\si{\sigma}

\def\CC{{\Bbb C}}

\def\NN{{\Bbb N}}
\def\PP{{\Bbb P}}

\def\RR{{\Bbb R}}
\def\TT{{\Bbb T}}
\def\ZZ{{\Bbb Z}}

\def\FSA{\Cal A}
\def\FSB{\Cal B}
\def\FSC{\Cal C}

\def\FSH{\Cal H}

\def\FSU{\Cal U}

\def\FSZ{\Cal Z}

\def\goh{{\goth h}}
\def\gog{{\goth g}}
\def\gok{{\goth k}}
\def\gol{{\goth l}}

\def\goo{{\goth o}}

\def\gos{{\goth s}}

\def\goS{{\goth S}}

\def\id{\operatorname{id}}
\def\phi{\varphi}

\def\qq{(q-q^{-1})}
\def\Aq{\FSA_q}
\def\Uq{\FSU_q}
\def\Bq{\FSB_q}
\def\detq{\det\nolimits_q}
\def\ks{{\gok^\sigma}}
\def\kss{(\gok^\sigma)^\ast}
\def\Lp{L^+}
\def\Lm{L^-}
\def\Lpm{L^\pm}
\def\Hst{\FSH^{(\si,\tau)}}
\def\xst{x^{(\si,\tau)}}
\def\AT{\FSA(\TT)}
\def\Uh{\Uq(\goh)}

\def\End{\operatorname{End}}
\def\Hom{\operatorname{Hom}}

\def\mez{{1\over 2}}

\def\innp#1 #2 {\langle #1 , #2 \rangle}
\let\eps=\varepsilon

\let\ten=\otimes

\documentstyle{amsppt} 
\def\AApart{AA1}
\def\AAclass{AA2}
\def\AW{AW}
\def\DKqsf{DK1}
\def\DKcqg{DK2}
\def\DZ{Dz}
\def\Drin{Dr}
\def\GR{GR}
\def\Jim{J}
\def\KWqjac{K1}
\def\Kor{K2}
\def\Kaw{K3}
\def\Koe{Kk}
\def\Lz{Lz}
\def\Mas{Ma}
\def\Noumimac{N}
\def\NMbigJ{NM1}
\def\NMqthreesphere{NM3}
\def\NMaw{NM2}
\def\NSug{NS1}
\def\NSb{NS2}
\def\NYM{NYM}
\def\Pod{P}
\def\RTF{RTF}
\def\Rosso{Ro}
\def\Sug{S1}
\def\Sb{S2}
\def\VKproj{KV}
\def\VS{VS}
\def\Wor{Wz}
\topmatter
\title   A family of quantum projective spaces \\
         and related  $q$-hypergeometric \\
         orthogonal polynomials
\endtitle
\author 
Mathijs S.\ Dijkhuizen and Masatoshi Noumi 
\endauthor
\affil Department of Mathematics, Faculty of Science,\\
Kobe University, Rokko, Kobe 657, Japan\endaffil
\email msdz$\@$math.s.kobe-u.ac.jp and
noumi$\@$math.s.kobe-u.ac.jp\endemail
\thanks The first author acknowledges financial 
support by the Japan Society for the Promotion
of Science (JSPS) and the Netherlands Organization 
for Scientific Research (NWO). \endthanks
\keywords quantum unitary group, quantum projective space,
two-sided coideal, zonal spherical function, Casimir operator, radial part,
second-order $q$-difference operator, Askey-Wilson polynomials,
big and little $q$-Jacobi polynomials \endkeywords
\subjclass 33D80, 81R50, 17B37, 33D45 \endsubjclass
\abstract We define a one-parameter family of 
two-sided coideals in $\Uq(\gog\gol(n))$ and 
study the corresponding algebras of 
infinitesimally right invariant functions on the 
quantum unitary group $U_q(n)$. The Plancherel 
decomposition of these algebras with respect to the
natural transitive $U_q(n)$-action is shown to be the 
same as in the case of a complex projective space. 
By computing the radial part of a suitable
Casimir operator, we identify the zonal spherical 
functions (i.e.\ infinitesimally bi-invariant 
matrix coefficients of finite-dimensional 
irreducible representations)
as Askey-Wilson polynomials containing two continuous
and one discrete parameter. In certain limit cases, the 
zonal spherical functions are expressed as big and 
little $q$-Jacobi polynomials depending on one 
discrete parameter.  
\endabstract
\rightheadtext{Quantum projective spaces}
\leftheadtext{Quantum projective spaces}
\endtopmatter
\document
\beginsektion 0 {Introduction}
In this paper, we study a family of two-sided coideals 
$\gok^{(c,d)}$ ($c,d$ non-negative real numbers) in the
quantized universal enveloping algebra $\Uq(\gog\gol(n))$.
The coideals $\gok^{(c,d)}$ can be viewed as a $q$-analogue 
of the Lie subalgebra 
$\gog\gol(n-1) \oplus \gog\gol(1) \subset \gog\gol(n)$.
By considering the algebra of functions on the quantum 
unitary group $U_q(n)$ that are ``infinitesimally'' 
invariant with respect to the
coideal $\gok^{(c,d)}$, we obtain a family of quantum 
projective spaces $\CC\PP^{n-1}_q(c,d)$ endowed with a 
natural transitive action of
the quantum unitary group $U_q(n)$.  These quantum 
$U_q(n)$-spaces were studied for the first time by 
Vaksman and Korogodsky \cite{\VKproj}
who defined 
them in a ``global'' way by means of a $q$-analogue 
of the classical Hopf fibration
$S^{2n-1} \to \CC\PP^{n-1}$. We also analyse the zonal 
spherical functions (infinitesimally left 
$\gok^{(c,d)}$-invariant and right 
$\gok^{(c',d')}$-invariant matrix coefficients) 
corresponding to finite-dimensional irreducible 
representations of $U_q(n)$.
They are expressed in terms of a family of 
Askey-Wilson polynomials
containing two continuous and one discrete parameter. 
We obtain this result
by showing that the zonal spherical functions are 
eigenfunctions of
a certain second-order $q$-difference operator 
which arises as the
radial part of a suitable Casimir operator.

The method of constructing  quantum homogeneous spaces 
by using the notion of ``infinitesimal'' invariance 
with respect to the natural action of
the quantized universal enveloping algebra was introduced
 by Koornwinder \cite{\Kor}, \cite{\Kaw} in the case of the
 quantum $SU(2)$ group. He actually worked with so-called
 twisted primitive elements in the
quantized universal enveloping algebra $\Uq(\gos\gol(2))$.
For other applications of this idea to the $SU_q(2)$ case 
we refer to \cite{\NMaw}, \cite{\Koe}, \cite{\DKqsf}.
 By extending the notion of twisted primitive element to
 the more general one of two-sided
coideal, it was shown in subsequent papers (cf.\  
\cite{\Noumimac}, \cite{\Sug}, \cite{\NSug}, see 
also \cite{\DKqsf})
that the infinitesimal method could also be fruitfully 
applied in the case of higher-dimensional quantum groups.
It not only provides one with a more practical
way to construct examples of quantum homogeneous spaces 
(finding a suitable two-sided coideal is easier than 
constructing the algebra of functions on
a quantum homogeneous space by means of generators 
and relations), but it also allows one to study infinitesimally 
bi-invariant functions or (zonal) spherical functions, 
something which can be done to only a very
limited extent if the quantum homogeneous space has been 
defined in a global way.

The results of this paper constitute a simultaneous 
generalization of results in various other papers. 
First of all, the case $c=1$, $d=0$ ($n$ arbitrary) 
was dealt with by Mimachi, Yamada and the second author 
\cite{\NYM}. In this case, one can
do everything globally, since the notions of infinitesimal
 invariance with respect to the coideal $\gok^{(1,0)}$ and
 global invariance with respect to the quantum subgroup 
$U_q(n-1)\times U(1)$ coincide.
The resulting spherical functions are 
little $q$-Jacobi polynomials, which are orthogonal 
 with respect to a measure supported on an 
infinite discrete set (Jackson $q$-integral).  

Secondly, the  case $n=2$ ($c,d$ arbitrary) 
has been thoroughly studied in a number
of papers. The quantum projective spaces then 
reduce to quantum spheres, which were introduced 
globally by Podle\'s \cite{\Pod}. The zonal spherical
functions on these quantum spheres were studied by 
Mimachi and the second author \cite{\NMbigJ}, and, 
in general, by Koornwinder \cite{\Kor}, \cite{\Kaw}, 
who used the infinitesimal method. For the special
case $n=2$ ($c=1, d=0$) we also refer to the papers
by Vaksman and So\u\i bel'man \cite{\VS}, Masuda et
al.\ \cite{\Mas}, and Koornwinder \cite{\KWqjac}.

The organization of this paper is as follows. In
 sections 1 and 2 we collect the necessary facts about the
 quantum unitary group and
$q$-hypergeometric orthogonal polynomials respectively. 
In section 3 we introduce the coideals $\gok^{(c,d)}$ and
 study the corresponding quantized
algebras of invariant functions and their Plancherel
 decomposition under the natural transitive $U_q(n)$-action. 
In section 4 we show that the algebras of invariant functions
can be constructed by means of a $q$-analogue of the
Hopf fibration $S^{2n-1} \to \CC\PP^{n-1}$ and establish
the link with Vaksman and Korogodsky \cite{KV}.
In section 5 we consider the zonal spherical
functions and identify them as a (partially discrete)
 three-parameter family of Askey-Wilson polynomials. 
We also treat some special cases ($c=0$ or
$d=0$) in which the Askey-Wilson polynomials degenerate 
to big or little $q$-Jacobi polynomials. In section 6
we present the details of the computation of the radial
part of a suitable Casimir operator.

The main results of this paper were announced in
\cite{\DZ}. The reader is referred to this last paper
for more background information and 
a discussion of the results presented here in
the context of general higher-rank quantum symmetric
spaces \cite{\Noumimac}, \cite{\NSug}, \cite{\NSb}, \cite{\Sb}.

The authors would like to thank T.\  Sugitani for 
numerous stimulating discussions and some very useful 
suggestions. The research for this paper
was started at the University of Amsterdam. 
The authors would like to express their gratitude to 
Prof.\ Tom H.\ Koornwinder for his hospitality and
for providing a stimulating environment in which to
 do research, the first author as a former temporary 
member of Prof.\ Koornwinder's research group, the
second one as an invited speaker at the Thomas
Stieltjes Research Institute during the Concentration Period on
Representation Theory and $q$-Special Functions (April-May 1994).
The first author would also like to thank the  
Mittag-Leffler Institute in Stockholm for its hospitality while 
preparing the final version of this paper.
\beginsektion 1 {Preliminaries on the quantum unitary group} 
The quantum unitary group $U_q(n)$ or its ``infinitesimal'' version,
the quantized universal enveloping algebra $\Uq = \Uq(\gog\gol(n))$,
have been studied in many papers, for instance
\cite{\Jim}, \cite{\Drin}, \cite{\RTF}, \cite{\NYM}, \cite{\Noumimac}.
Our basic reference will be \cite{\Noumimac, \S 1}.
In this section we only recall the barest essentials. Our notation
is virtually the same as in \cite{\Noumimac, \S 1}.

Let us fix $0<q<1$ and $n\geq 2$. Let $V$ be the vector space over $\CC$ 
with canonical basis $(v_i)_{1\leq i\leq n}$. 
Our starting point is the  
invertible $n^2\times n^2$ matrix $R\in\End(V\ten V)$ defined by
$$R := \sum_{ij} q^{\de_{ij}} e_{ii}\ten e_{jj} +
\qq \sum_{i>j} e_{ij}\ten e_{ji},\eqno\eq{Rdef}$$
where the $e_{ij}\in\End(V)$ denote the standard 
matrix units with respect to the basis $(v_i)$. 
The generators $t_{ij}$ of $\Aq = \Aq(U(n))$ satisfy the usual
commutation relations $R T_1 T_2 = T_2 T_1 R$ and form a
unitary matrix corepresentation of $\Aq$.

Let $\AT:=\CC[z_1^{\pm 1}, \ldots, z_n^{\pm 1}]$
be the algebra of trigonometric polynomials on the
$n$-dimensional real torus $\TT$, the Hopf $\ast$-algebra structure
being given by  
$$\Delta(z_i) = z_i \ten z_i, \quad \eps(z_i) = 1, \quad
z_i^\ast = z_i^{-1} \quad (1\leq i\leq n).
\eqno\eq{comulttorus}$$
$\TT$ is naturally identified with the diagonal
subgroup of $U_q(n)$,
the corresponding surjective Hopf $\ast$-algebra morphism
(restriction of functions) being denoted by
$${}_{|\TT}\colon \Aq \longrightarrow \AT.\eqno\eq{torus-restrict}$$

Let $P = \bigoplus_{1\leq i\leq n} \ZZ\eps_i$ denote the weight lattice
of $U(n)$. We identify $P$ and $P^\ast = \Hom_\ZZ(P,\ZZ)$ by means 
of the pairing $\langle \eps_i, \eps_j \rangle = \de_{ij}$. 
The algebra $\Uq$ is generated by the symbols 
$q^h$ ($h\in P^\ast$) and $e_i, f_i$ ($1\leq i \leq n-1$) 
subject to the well-known quantized Weyl-Serre relations.
Let $\Uh\subset \Uq$ denote the subalgebra generated by
the elements $q^h$ ($h\in P^\ast$). 

We put
$$R^+ := PRP, \qquad R^- := R^{-1}, \eqno\eq{Rplusmin}$$
where $P\in\End(V\ten V)$ is the usual permutation operator.
One has the identities $(R^\epsilon)^{-1} = 
(R^{-\epsilon})^t$ ($\epsilon = \pm$).
The algebra $\Uq$ is also generated by 
the so-called {\sl $L$-operators}
 $\Lp_{ij},\,\Lm_{ij}\in \Uq$ (cf.\ \cite{\Jim}, \cite{\RTF},
\cite{\Noumimac}) subject to the relations 
$$R^+ L_1^\epsilon L_2^\epsilon = L_2^\epsilon L_1^\epsilon R^+\; 
(\epsilon = \pm), \quad 
R^+ \Lp_1 \Lm_2 = \Lm_2 \Lp_1  R^+, \eqno\eq{relL}$$
where $L^\pm := (L^\pm_{ij})$, $\Lpm_1 := \Lpm \ten \id$ etc.
In particular,
$$ q^h \Lpm_{ij} q^{-h} = q^{\langle  h,  \eps_j -\eps_i\rangle} 
\Lpm_{ij}, \quad
q^h S(\Lpm_{ij}) q^{-h} = q^{\langle h, \eps_j -\eps_i\rangle}
S(\Lpm_{ij})\; (1\leq i,j\leq n). \eqno\eq{qhLrel}$$
The $L^{\pm}_{ij}$ can be viewed
as quantum analogues of arbitrary root vectors in $\gog\gol(n)$.
The matrices $L^+$ and $L^-$ are upper and lower triangular
respectively.

The Hopf $\ast$-algebra structure on
$\Uq$ is uniquely determined by 
$$\Delta(L^\pm_{ij}) = \sum_k L^\pm_{ik} \ten 
L^\pm_{kj}, \quad \eps(L^\pm_{ij}) = \delta_{ij},
\quad (L^\pm_{ij})^\ast = S(L^\mp_{ji}) \quad (1\leq i,j\leq n)
\eqno\eq{deltaL}$$

The action of the involution
$\tau = \ast \circ S\colon \Uq \to \Uq$ on the generators is given by
$$\tau(\Lpm_{ij}) = L^\mp_{ji}\quad
(1\leq i,j \leq n).\eqno\eq{tauL}$$

The cone $P^+\subset P$ of {\sl dominant weights} by definition
consists of all weights $\lambda = \sum_k \lambda_k\eps_k \in P$ such that
$\lambda_1 \geq \ldots \geq \lambda_n$.
There is the usual parametrization $\lambda \mapsto V(\lambda)$
of irreducible $P$-weighted finite-dimensional left $\Uq$-modules
by dominant weights (cf.\ \cite{\Lz}, \cite{\Rosso}). Recall that 
a highest weight vector  of a left $\Uq$-module $W$
is  by definition annihilated by the $L^-_{ij}$ or, equivalently,
by the $S(L^-_{ij})$ ($i>j$).
All finite-dimensional $P$-weighted left $\Uq$-modules are completely 
reducible and unitarizable.

Weights of right $\Uq$-modules are defined in the obvious way, 
highest weight vectors being killed by definition by the 
$L^+_{ij}$ or the $S(L^+_{ij})$ ($i>j$).
Given a left $\Uq$-module $W$, one defines a right $\Uq$-module
structure on the conjugate vector space $W^\circ$ by 
putting
$$ v \cdot u := u^\ast \cdot v 
\quad (v\in W, u\in \Uq).\eqno\eq{leftright}$$
The assignment $W\mapsto W^\circ$ is a 1-1 correspondence
between left and right $\Uq$-modules 
preserving weight vectors, weights,  and highest weights.

There exists a unique algebra homomorphism $\rho_V\colon
\Uq \to \End(V)$ such that
$$R^\pm = \sum_{ij} e_{ij} \ten \rho_V(\Lpm_{ij}), \quad
(R^\pm)^{-1} = \sum_{ij} e_{ij} \ten \rho_V(S(L^\pm_{ij})).
\eqno\eq{vrepdef}$$ 
The corresponding  representation 
$V$ is called {\sl vector representation} and has 
highest weight $\eps_1$. The elements $q^h\in \Uq(\goh)$
act on $V$ as diagonal matrices with respect to the basis
$(v_i)$.

By means of the natural Hopf $\ast$-algebra 
duality $\langle \cdot\, , \, \cdot \rangle$
 between $\Uq$ and $\Aq$, one can identify $\Aq$ with a 
subspace of the algebraic linear dual of $\Uq$. 
There is an induced Hopf $\ast$-algebra duality between
$\Uq(\goh)$ and $\AT$ such that 
$$\langle q^h\,,\, z^\lambda\rangle 
 := q^{\langle h, \lambda\rangle}, \quad
z^\lambda = z_1^{\lambda_1} \cdots z_n^{\lambda_n}\quad
(h\in P^\ast,\; \lambda =\sum_k \lambda_k\eps_k \in P).
\eqno\eq{dualtorus}$$

The duality naturally turns $\Aq$ into a $\Uq$-bimodule with
two-sided $\Uq$-symmetry.
Recall that transposition defines a natural right $\Uq$-module
structure on the linear dual $\Hom(V,\CC)=V^\ast$. 
The mapping 
$$V\ten \Hom(V,\CC) \to \Aq, \quad 
v_i \ten v_j^\ast  \mapsto t_{ji} \eqno\eq{Wvector}$$
is an (injective) $\Uq$-bimodule homomorphism.
More generally, one has the following 
decomposition of $\Aq$ into irreducible $\Uq$-bimodules:
$$\Aq = \bigoplus_{\lambda\in P^+} 
V(\lambda) \ten \Hom(V(\lambda),\CC).
\eqno\eq{Aqdecomp}$$
Here the subspace 
$W(\lambda) := V(\lambda) \ten \Hom(V(\lambda),\CC)
\subset \Aq$
is spanned by the matrix coefficients of the (co-)representation
$V(\lambda)$. Note that $\Hom(V(\lambda), \CC)$ is isomorphic
with $V(\lambda)^\circ$, although the isomorphism is not canonical.
The decomposition \eqtag{Aqdecomp} can also be characterized as the
simultaneous eigen\-space decomposition of $\Aq$ with respect to the
natural action of the center $\FSZ\Uq\subset \Uq$.

Let $h\colon \Aq\to\CC$ denote the Haar functional on $\Aq$
(cf.\ \cite{\Wor}, \cite{NYM}, \cite{\DKcqg}). 
Then $\langle a, b\rangle :=
h(b^\ast a)$ defines a positive definite inner product on $\Aq$ 
with respect to which the subspaces $W(\lambda)\subset \Aq$ are
mutually orthogonal (Schur orthogonality). 
\beginsektion 2 {Preliminaries on $q$-hypergeometric orthogonal
polynomials}
Let $0<q<1$. The so-called {\sl $q$-shifted factorials}
are defined as
$$(a;q)_n := \prod_{k=0}^{n-1} (1-aq^k), \;
(a_1,\ldots, a_s;q)_n := \prod_{j=1}^s (a_j;q)_n,
\; (a;q)_\infty := \lim_{n\to\infty} (a;q)_n,\eqno\eq{qfactorial}$$
and the {\sl $q$-hypergeometric series} 
${}_{s+1}\phi_s$ as
$${}_{s+1}\phi_s \left [ {a_1, \ldots, a_{s+1} \atop
b_1, \ldots, b_s}; q,z \right ] :=
\sum_{k=0}^\infty {(a_1, \dots, a_{s+1};q)_k \,z^k\over
(b_1, \ldots, b_s;q)_k\, (q;q)_k}.\eqno\eq{qhypseries}$$
We shall only deal with the special case of \eqtag{qhypseries}
when $a_1 = q^{-n}$ ($n\in \ZZ_+$). The series 
then is {\sl terminating}, i.e.\ for $j>n$ the $j$-th
term vanishes. In this case, it is always tacitly assumed that
$b_1, \ldots, b_s \notin \{1,q^{-1}, \ldots, q^{-n+1}\}$.
Further details about $q$-hypergeometric functions can be found
in \cite{\GR, Ch.\ 1}.

Consider the Laurent polynomial ring $\CC[z^{\pm 1}]$
in the variable $z$ and put $x:= {1\over 2} (z+ z^{-1})$.
The {\sl Askey-Wilson polynomials} (cf.\ \cite{\AW}) 
are the polynomials in the variable $x$ given by
$$\eqalignno{p_n(x &; a,b,c,d \mid q) := & \cr
&a^{-n} (ab,ac,ad;q)_n  \cdot
{}_4\phi_3\left [ {q^{-n}, q^{n-1}abcd, az, 
az^{-1}\atop
ab,ac,ad}; q,q \right ].& \eq{defaw}\cr}$$
They are symmetric in $a,b,c,d$. It is sometimes convenient
to write the ${}_4\phi_3$ factor in \eqtag{defaw} as
$r_n(x) = R_n(z)$.
Depending on the value of
$a,b,c,d$, Askey-Wilson polynomials are orthogonal
 with respect to a positive  orthogonality measure consisting of 
a continuous and a discrete part.  To be more precise,
one has the following result (cf.\ \cite{\AW, Thm. 2.5}):
\beginth{Proposition} orthaw
{\thstil Assume that (i) $a,b,c,d$ are real, or, if complex, appear
in conjugate pairs, and  (ii) the pairwise products of
$a,b,c,d$ are not $\geq 1$. Then the Askey-Wilson polynomials 
$p_n(x)$ are mutually orthogonal with respect to the positive
definite inner product on $\CC[x]$ defined by
$$\eqalignno{\langle P, Q \rangle :=
{1\over 2\pi i} \int_{z\in\FSC} P(x)\, &\overline{Q(x)} \,w(z; a,b,c,d;q)
\,{\hbox{d} z \over z} & \eq{orthaw-rel} \cr
& + \sum_{e,k} P(e_k) \overline{Q(e_k)}w_k(e,f,g,h)\quad
(P,Q\in \CC[x]).& \cr}$$
Here integration is along the unit circle $\FSC$ in counterclockwise
direction, the summation runs over all $e\in\{a,b,c,d\}$ such that $|e|>1$
and all $k\in \ZZ_+$ such that $|eq^k| > 1$. Moreover, the following
notation is used:
$$w(z; a,b,c,d;q) := {(z^2, z^{-2};q)_\infty\over
(az, a/z, bz, b/z, cz, c/z, dz, d/z;q)_\infty},$$
$$\eqalign{w_k(e,f,g,h) := 
&{(e^{-2};q)_\infty \over (q,ef,f/e,eg,g/e,eh,h/e;q)_\infty}\times \cr
& \qquad {(e^2,ef,eg,eh;q)_k\over (q,eq/f,eq/g,eq/h;q)_k}
{(1-e^2q^{2k}) \over(1-e^2)} \left ( {q\over efgh} \right )^k,\cr}$$
%$$e_k :=
%(eq^k + e^{-1}q^{-k})/2, \quad \{f,g,h\} := \{a,b,c,d\}\backslash \{e\}.$$
$$e_k := (eq^k + e^{-1}q^{-k})/2,$$
and $(e,f,g,h)$ is any permutation of $(a,b,c,d)$.
\par}
\endth\noindent
One can calculate an explicit expression for 
the norm of the Askey-Wilson
polynomials $p_n(x; a,b,c,d \mid q)$ with 
respect to the inner product
\eqtag{orthaw-rel} (cf.\ \cite{\AW, Thm. 2.5}). 
We shall only need the value for $p_0 = 1$ (cf.\ 
\cite{\AW, Thm.\ 2.1}):
$$\langle 1, 1 \rangle = {2(abcd;q)_\infty \over
(q, ab, ac, ad, bc, bd, cd; q)_\infty}.
\eqno\eq{awnorm}$$
 For 
certain special values of the parameters, e.g.\ 
$\max(|a|, |b|, |c|, |d|) < 1$, 
the discrete part in \eqtag{orthaw-rel}
becomes void. The orthogonality measure 
then is absolutely continuous.

The Laurent polynomials $R_n(z)$ in the variable $z$
satisfy the following second-order $q$-difference
equation:
$$\eqalignno{A&(z;q)(R_n(qz) - R_n(z)) + 
A(z^{-1};q)(R_n(q^{-1}z) -R_n(z)) = & \eq{qdiffeq}\cr
&= -(1-q^{-n})(1-q^{n-1}abcd)R_n(z), & \cr}$$
where
$$A(z;q) := 
{(1-az)(1-bz)(1-cz)(1-dz)\over (1-z^2)(1-qz^2)} \quad
(a,b,c,d\in\CC).\eqno\eq{Azqdef}$$
Any symmetric Laurent polynomial $f(z)$ that is of degree $\leq n$ 
when viewed as a polynomial in $x$  and satisfies
 \eqtag{qdiffeq}  is a constant multiple of $R_n(z)$.

There are the so-called 
{\sl big $q$-Jacobi
polynomials} (cf.\ \cite{\AAclass}):
$$P_n^{(\al,\be)}(x;c,d\colon q) :=
{}_3\phi_2\left [ {q^{-n}, q^{n+\al+\be+1}, q^{\al+1}x/c \atop
q^{\al+1}, -q^{\al+1}d/c}; q,q \right ]. \eqno\eq{defbigqjac}$$
If $c,d>0$ and $\alpha, \beta > -1$, 
then the polynomials $P_n^{(\al,\be)}(x;c,d\colon q)$ are orthogonal
with respect to a positive orthogonality measure supported on the
infinite discrete set 
$$\{cq^k \mid k\in \ZZ_+\} \cup \{-dq^k \mid k\in \ZZ_+\} \subset
[-d,c].$$

A special case of the big $q$-Jacobi polynomials ($c=1,d=0$)
are the so-called {\sl little
$q$-Jacobi polynomials} (cf.\ \cite{\AApart}):
$$p_n^{(\al,\be)}(x\colon q) := {}_2\phi_1(q^{-n}, 
q^{\al+\be + n+1}; q^{\al + 1}; q; qx). 
\eqno\eq{deflittleqjac}$$
For $\al, \be > -1$ they are orthogonal 
with respect to a positive orthogonality measure supported on
the infinite discrete set $\{q^k \mid k\in \ZZ_+ \} \subset [0,1]$.

Both big and little $q$-Jacobi polynomials 
may be characterized as the polynomial
eigenfunctions of a certain second-order $q$-difference operator
with rational coefficients and depending on the parameters
$\alpha$, $\beta$ (and $c,d$ in the big $q$-Jacobi case).

Recall the notation $r_n(x)= r_n(x;a,b,c,d\mid q)$ for the 
${}_4\phi_3$ factor in \eqtag{defaw}.
Little and big $q$-Jacobi polynomials 
can be recovered from the Askey-Wilson polynomials $r_n$ by a suitable
limit transition (cf.\ \cite{\Kaw, Prop.\  6.1, 6.3}):
\beginth{Proposition} limitbig
{\thstil Let the big 
$q$-Jacobi polynomials be denoted as in 
\eqtag{defbigqjac}.
Then 
$$\eqalign{\lim_{a\to 0}\,r_n \biggl ( &{q^{\mez} x\over 
2a(cd)^{\mez}}; q^{\al+{1\over 2}} a (d/c)^{\mez}, 
q^{\mez} a^{-1} (c/d)^{\mez},\cr
&-q^{\mez} a^{-1} (d/c)^{\mez}, 
-q^{\be+{\mez}} a (c/d)^{\mez} 
\mid q \biggr ) = P_n^{(\al,\be)}(x; c,d\colon q).\cr}$$
\par}
\endth
\beginth{Proposition} limitlittle
{\thstil Let the little 
$q$-Jacobi polynomials be denoted as in 
\eqtag{deflittleqjac}.
Then 
$$\lim_{a\to 0}\,r_n \biggl ({q^{\mez} x\over 2a^2};
\pm q^{\al+\mez} a^2, \pm q^\mez a^{-2}, \mp q^\mez, 
\mp q^{\be+\mez} \mid q \biggr ) =
{(q^{\be+1};q)_n \over (q^{-n-\al};q)_n} 
p_n^{(\be,\al)}(x\colon q).$$
\par}
\endth\noindent
As follows from Proposition \thtag{orthaw} after a suitable transformation
of variables (cf.\ \cite{\Kaw, Remark 6.6}), 
for $\al, \be >-1$, $c,d >0$ and $a$ small enough (and real), the orthogonal
polynomials in $x$ on the left-hand side of the displayed
equation in Proposition \thtag{limitbig} (after the limit sign) have
continuous mass on the interval $[-2a (cd/q)^{1\over 2},
2a(cd/q)^{1\over 2}]$ and discrete mass on the two
finite sets
$$\{cq^k + a^2 d q^{-k-1} \mid k\in \ZZ_+\,,\, 
q^k > a(qc/d)^{-{1\over 2}}\}$$
and 
$$\{-dq^k - a^2 c q^{-k-1} \mid k\in \ZZ_+\,,\, 
q^k > a(qd/c)^{-{1\over 2}}\}.$$
When $a\to 0$, the continuous mass interval shrinks to
$\{0\}$, while the two discrete sets tend to the support
of the orthogonality measure of the big $q$-Jacobi polynomials.
A similar remark applies to little $q$-Jacobi polynomials.
\beginsektion 3 {A family of two-sided coideals}
Let us fix real numbers $c,d\geq 0$ such that 
$(c,d) \neq (0,0)$. 
The subspace 
$\gok^{(c,d)} \subset \Uq$ is by 
definition spanned 
by the following elements:
$$\eqalignno{(i) \quad & \Lp_{11} - \Lm_{nn}, \;
\Lm_{11} - \Lp_{nn}, & \cr
(ii) \quad & \sqrt{c}\, \Lp_{1k} + \sqrt{d}\, 
\Lm_{nk} \quad (2\leq k\leq n-1), & \cr
(iii) \quad & \sqrt{d}\, \Lp_{kn} + 
\sqrt{c}\, \Lm_{k1} \quad
(2\leq k\leq n-1), & \eq{kdefcd} \cr
(iv) \quad & \Lp_{ij}, \; \Lm_{ji} \quad 
(2\leq i<j\leq n-1), & \cr
(v) \quad & \Lp_{ii} - \Lm_{ii} \quad 
(2\leq i \leq n-1), & \cr
(vi) \quad & \sqrt{cd}\, \Lp_{1n} - 
\sqrt{cd}\, \Lm_{n1}
-(c-d) (\Lp_{11} -\Lm_{11}). & \cr}$$
We remark that the subspace $\gok\subset \Uq$ only 
depends on the ratio of the numbers $c$ and $d$.
In fact, it will be convenient to introduce a parameter
$\sigma\in \RR\cup \{\pm \infty\}$ by setting
$$q^\sigma = \sqrt{{d\over c}} \;\; (c,d>0),\quad
\sigma = \infty\;\; (d=0),\quad \sigma = -\infty\;\; (c=0).
\eqno\eq{sigmadef}$$
Then $\gok^{(c,d)}$ only depends on $\sigma$ and 
we write $\gok^\sigma := \gok^{(c,d)}$.

In the remainder of this paper, whenever we write
$\gok^\sigma$ or $\gok^{(c,d)}$, it is tacitly
assumed that $-\infty \leq \sigma \leq \infty$
and $c,d \geq 0$, $(c,d)\neq (0,0)$, unless 
explicitly mentioned otherwise.
\beginth{Proposition} kcoideal
{\thstil 
The subspace $\gok^\sigma\subset \Uq$ is a 
$\tau$-invariant two-sided coideal.
\par}
\endth
\beginproof
It follows from \eqtag{tauL} that $\tau$ permutes the elements
listed in \eqtag{kdefcd} (up to a scalar multiple). Hence
$\gok^\sigma$ is $\tau$-invariant.
The property $\eps(\ks)= 0$ is a direct consequence of
\eqtag{deltaL}. A straightforward computation using 
\eqtag{deltaL} shows that $\Delta(\ks) \subset
\ks\ten \Uq + \Uq \ten \ks$. We leave the details to
the reader.
\endproof
We now proceed to study $\ks$-fixed and $\kss$-fixed
vectors in finite-dimensional
representations of $\Uq$.
Let $W$ be a left $\Uq$-module. A vector $v\in W$ is
called {\sl invariant} 
w.r.t.\ an element
$u\in\Uq$ if $u\cdot v = \eps(u) v$. In particular, the subspace
 $W_\ks \subset W$ of {\sl $\ks$-fixed vectors} 
is defined as 
$$W_\ks := \{ v\in W \mid \ks \cdot v = 0 \}.\eqno\eq{kfixeddef}$$
It is obvious how to define $\ks$-fixed vectors for
right $\Uq$-modules. Note that a vector $v$ in a right
$\Uq$-module $W$ is $\ks$-fixed if and only if $v$ is
$(\ks)^\ast$-fixed as an element of the left $\Uq$-module
$W^\circ$. 

Assuming that $\sigma$ is
finite, we have the following crucial result:
\beginth{Proposition} hwcomp
{\thstil Suppose $\sigma$ is finite and let $\gok$ denote either
$\ks$ or $\kss$. Suppose $\lambda\in P^+$. If $v\in V(\lambda)$ is 
a non-zero $\gok$-fixed vector 
then the highest weight component of $v$ is non-zero.
\par}
\endth
\beginproof
The proof is similar to that of \cite{\Noumimac, Lemma 3.2}.
Suppose that $0\neq v\in V(\lambda)$ is $\ks$-fixed and write 
$v$ as a sum 
$v=\sum_{\mu\in P} v_\mu$ of weight vectors.
Recall (cf.\ \eqtag{qhLrel}) that $L^\pm_{ij}$ has
weight $\eps_j -\eps_i$ in $\Uq$. This implies that
$L^\pm_{ij}\cdot v_\mu\in V(\lambda)$ has weight $\mu+\eps_j-\eps_i$.
Let $\leq$ denote the lexicographic order on $P$ with respect
to the $\ZZ$-basis $(\eps_i)$. One reads off from \eqtag{kdefcd}
that, for any $i<j$,
there is an element $X_{ij}\in \ks$ whose leading term (w.r.t.\  $\leq$)
 is equal to $L^-_{ji}$. Let now
$\mu_0\in P$ be the greatest element $\mu \in P$ (w.r.t.\ $\leq$)
 such that $v_\mu \neq 0$. Then $L^-_{ji}\cdot v_{\mu_0}$
is the component of
weight $\mu_0 +\eps_j -\eps_i$ in $X_{ij}\cdot v$. But $X_{ij}
\cdot v =0$, since $X_{ij} \in \ks$. Hence $L^-_{ij}\cdot
v_{\mu_0} = 0$. In other words, $v_{\mu_0}$ is a non-zero
highest weight vector in $V(\lambda)$, which forces $\mu_0 = \lambda$.
One can prove the corresponding statement
for $\kss$ in a completely analogous way using \eqtag{deltaL} and
the fact that a highest weight vector $v$ is characterized among
weight vectors by the condition $S(L^-_{ij}) \cdot v = 0$ ($1\leq
j<i\leq n$).
\endproof
\beginth{Corollary} gelfand
{\thstil Suppose $\sigma$ is finite and let $\gok$ denote either
$\gok^\sigma$ or $\kss$.  Suppose $\lambda\in P^+$. The subspace
$V(\lambda)_\gok$ of $\gok$-fixed vectors is at most one-dimensional.
\par}
\endth
\beginproof
This is a direct consequence of the preceding proposition
(cf.\ \cite{\Noumimac, \S 3.1}).
\endproof
Corollary \thtag{gelfand} can also be expressed by saying
that $(\Uq, \ks)$ resp.\ $(\Uq, \kss)$ satisfy the 
Gelfand pair property. We call $V(\lambda)$ ($\lambda\in P^+$)
{\sl spherical} (with respect to $\gok^\sigma$)
if it has non-zero $\ks$-fixed vectors.

In order to determine the spherical representations
we need some explicit information about the vector
representation $V$ and its contragredient $V^\ast$.
One deduces from \eqtag{vrepdef} that the action of
the $L^\pm_{ij}$ and the $S(L^\pm_{ij})$ 
on the basis vectors $v_i\in V$ is given by the following explicit
formulae:
$$\eqalignno{L^\epsilon_{ii} \cdot v_k 
&= q^{\epsilon \de_{ik}} v_k, \quad
S(L^\epsilon_{ii}) \cdot v_k = q^{-\epsilon \delta_{ik}}v_k
\; (\epsilon = \pm), & \eq{Lvector}\cr
L^\pm_{ij}\cdot v_k &= \pm (q-q^{-1})\de_{ik}v_j, \quad
S(L^\pm_{ij})\cdot v_k = \mp(q-q^{-1})\de_{ik}v_j\; 
 (i \lessgtr j).& \cr}$$
The vector $v_i$ ($1\leq i\leq n$) has weight $\eps_i$. 
The highest weight vector is $v_1$.
Recall that the action of $\Uq$ on the contragredient module
$V^\ast$ is given by
$$(u\cdot v^\ast)(v) = v^\ast(S(u)\cdot v)\quad 
(v\in V, v^\ast\in V^\ast).\eqno\eq{Vdualdef}$$
Let $(v_i^\ast)$  denote the dual basis of $V^\ast$. The action
of the $L^\pm_{ij}$ and the $S(L^\pm_{ij})$ on the basis vectors
$v_k^\ast$ is given by:
$$\eqalignno{L^\epsilon_{ii} \cdot v_k^\ast &= 
q^{-\epsilon \de_{ik}} v_k^\ast,\quad
S(L^\epsilon_{ii}) \cdot v_k^\ast = 
q^{\epsilon \de_{ik}} v_k^\ast \; (\epsilon = \pm), & \eq{Lvectordual}\cr
L^\pm_{ij}\cdot v_k^\ast &= \mp(q-q^{-1})\de_{jk}v_i^\ast,
\; S(L^\pm_{ij})\cdot v_k^\ast = \pm q^{2(j-i)}(q-q^{-1})\de_{jk}
v_i^\ast\;(i \lessgtr j). & \cr}$$
To compute the action of the $S(L^\pm_{ij})$ in \eqtag{Lvectordual}
one needs the fact (cf.\ \cite{\Noumimac, (1.15)}) that the
square of the antipode $S^2\colon \Uq \to \Uq$ is
given by
$$S^2(u) = q^{-2\rho} u q^{2\rho}, \quad
\rho := \sum_{k=1}^n (n-k) \eps_k \in P. \eqno\eq{Ssquare}$$
The vector $v_i^\ast$ ($1\leq i\leq n$) has weight $-\eps_i$. 
The highest weight vector is $v_n^\ast$.
The $t_{ij}^\ast = S(t_{ji})$ are the coefficients of $V^\ast$
with respect to the basis $(v_i^\ast)$.

The tensor product $V^\ast\ten V$ has the following irreducible
decomposition:
$$V^\ast \ten V \cong V(0) \oplus V(\eps_1 - \eps_n).\eqno
\eq{VVstar}$$ 
Here  the subspace $V(0)$ is spanned
by the element $\sum_k q^{2(n-k)}v_k^\ast \ten v_k$.
The subspace $V(\eps_1-\eps_n)$ is spanned by the
linearly independent vectors
$$v_i^\ast \ten v_j \quad (1\leq i\neq j\leq n), \qquad 
v_i^\ast \ten v_i -v_{i+1}^\ast \ten v_{i+1} \quad
(1\leq i\leq n-1).\eqno\eq{Vepsbasis}$$
The vector $v_i^\ast\ten v_j$ ($1\leq i,j\leq n$) has
weight $\eps_j-\eps_i$. The highest weight vector
in $V(\eps_1\ten \eps_n)$ is $v_n^\ast
\ten v_1$.
All these statements can be easily deduced from 
\eqtag{Lvector}, \eqtag{Lvectordual}, and \eqtag{deltaL}.
\beginth{Theorem} spherrep
{\thstil Let $\sigma$ be finite and let $\gok$ denote either
$\gok^\sigma$ or $\kss$.  For any $\lambda\in P^+$,
the representation $V(\lambda)$ has non-zero $\gok$-fixed
vectors if and only if 
$\lambda = l(\eps_1 - \eps_n)$ for some $l\in \ZZ_+$.
\par}
\endth
\beginproof
 We first prove
the ``only if'' part of the statement. Let $v\in V(\lambda)$
be a non-zero $\gok$-fixed vector with highest weight component
$v_\lambda\neq 0$. For any $2\leq i\leq n-1$ one has
$0=(L^+_{ii} - L^-_{ii}) \cdot v_\lambda = 
(q^{\langle \lambda, \eps_i\rangle} - q^{-\langle\lambda,
\eps_i\rangle}) v_\lambda.$
Hence $\lambda \in \ZZ\eps_1 \oplus \ZZ\eps_n$. On the other
hand, $(L^+_{11}-L^-_{nn}) \cdot v_\lambda = 0$. This implies
$q^{\langle \lambda, \eps_1\rangle} =q^{-\langle\lambda,
\eps_n\rangle}$, hence $\lambda = l(\eps_1-\eps_n)$ for some
$l\in \ZZ_+$. 
To prove the ``if'' part of the statement we first exhibit
a $\gok$-fixed vector in the ``lowest'' spherical
representation $V(\eps_1 -\eps_n)$. In the representation
$V^\ast \ten V$ we have the following $\gok$-fixed vector
$w_\gok$:
$$\eqalignno{w_{\gok^\sigma} & 
:= \sqrt{cd}\,v^\ast_1 \ten v_n + \sqrt{cd}\, 
v^\ast_n \ten v_1 + q d\, v^\ast_1 \ten v_1 + 
q^{-1} c\, v^\ast_n \ten v_n, & \eq{repkfixed}\cr
w_{\kss} & := q^{2(n-1)}\sqrt{cd}\,v^\ast_1 \ten v_n + 
\sqrt{cd}\, v^\ast_n \ten v_1 + 
q^{2(n-1)}d\, v^\ast_1 \ten v_1 +  c\, v^\ast_n \ten v_n. & \cr}$$
One verifies that $w_\gok$ is $\gok$-fixed by means of
a straightforward computation using \eqtag{Lvector} and
\eqtag{Lvectordual}. Since the highest weight term in
$v_\gok$ is $\sqrt{cd}\, v_n^\ast \ten v_1
\neq 0$, the component of $v_\gok$ in $V(\eps_1-\eps_n)$
is non-zero and obviously $\gok$-fixed. To prove the
``if'' part in general, let us first remark that
the vector $(v_n^\ast \ten v_1)^{\ten l} \in (V^\ast \ten V)^{\ten l}$
is a highest weight vector of weight $l(\eps_1-\eps_n)$
occurring with multiplicity 1 inside $(V^\ast \ten V)^{\ten l}$.
It generates a $\Uq$-submodule isomorphic with 
$V(l(\eps_1-\eps_n))$. The component of weight $l(\eps_1-\eps_n)$
in $v_\gok^{\ten l}$ is equal to $(cd)^{l/2} 
(v_n^\ast \ten v_1)^{\ten l}$ which is obviously non-zero. 
This proves that there exists a non-zero $\gok$-fixed vector
in $V(l(\eps_1-\eps_n))$.
\endproof
The analysis of the spherical representations for 
$\sigma = \pm \infty$ proceeds along rather different lines,
since there is no analogue for Proposition \thtag{hwcomp}.
There is a natural embedding of Lie algebras
$$\gog\gol(n-1) \oplus \gog\gol(1) \hookrightarrow \gog\gol(n),
\quad (X, \xi) \mapsto \pmatrix X & 0 \\ 0 & \xi\endpmatrix .
\eqno\eq{gl-embedding}$$
This mapping arises from a natural embedding of
the Dynkin diagram of root system $A_{n-2}$ into that
of $A_{n-1}$. Since the definition of $\Uq(\gog\gol(n))$
is entirely in terms of generators corresponding to simple
roots and relations involving only the Cartan integers, 
there is an analogue of \eqtag{gl-embedding} for
quantized universal enveloping algebras. To be more
precise, by abuse of notation 
let $U(\gog\gol(1))$ denote the commutative algebra 
generated by the symbols $\zeta^{\pm 1}$ subject to the
relation $\zeta \zeta^{-1} = 1$. There is a unique
Hopf $\ast$-algebra structure on $U(\gog\gol(1))$ such that
$$\Delta(\zeta^{\pm 1}) = \ze^{\pm 1}\ten \ze^{\pm 1}, 
\quad \eps(\ze^{\pm 1}) = 1, \quad (\ze^{\pm 1})^\ast = \ze^{\pm 1}.
\eqno\eq{glonedef}$$
There is a natural  injective Hopf $\ast$-algebra homomorphism
$$\Uq(\gok):=\Uq(\gog\gol(n-1)) \ten U(\gog\gol(1))
\hookrightarrow \Uq(\gog\gol(n)) \eqno\eq{q-embedding}$$
sending $1 \ten \ze^{\pm 1}$ to $q^{\pm \eps_n}$ and
$L^\pm_{ij}\ten 1$ to $L^\pm_{ij}$ ($1\leq i,j \leq n-1$). 
Henceforth, we identify $\Uq(\gok)$ with its image
under this mapping. It is straightforward to check that
the notions of invariance w.r.t.\ $\Uq(\gok)$ and
$\gok^\infty$  coincide. 
The same is true for $(\gok^\infty)^\ast$, since 
$\Uq(\gok)$ is obviously $\ast$-invariant.

Dual to \eqtag{q-embedding} there 
is a natural surjective Hopf $\ast$-algebra isomorphism
$$\Aq(U(n)) \twoheadrightarrow
\Aq(U(n-1)) \otimes \FSA(U(1)) =: \Aq(K).\eqno\eq{K-embedding}$$
Here $\FSA(U(1))$ is the algebra of trigonometric polynomials
on the one-dimensional real torus $U(1)$ (cf.\ \eqtag{comulttorus}).

The natural notion of invariance w.r.t.\ 
$\Aq(K)$ coincides with
$\Uq(\gok)$-invariance (cf.\ \cite{\DKqsf, Prop. 1.12})
and hence with $\gok^\infty$-invariance.
$\Aq(K)$-fixed vectors for finite-dimensional corepresentations
of $\Aq$ have essentially been analysed in \cite{\NYM,\S 4}. 
It can be shown that for finite $\sigma$ there is no quantum subgroup
corresponding to the coideal $\gok^\sigma$.
\beginth {Remark} Dynkinautom
There is a similar picture in  the
case $\sigma = -\infty$. Instead of \eqtag{gl-embedding}
one starts with the embedding
$$\gog\gol(1)\oplus \gog\gol(n-1)\hookrightarrow \gog\gol(n),
\quad (\xi, X) \mapsto \pmatrix \xi & 0 \\ 0 & X \endpmatrix .$$
This embedding is conjugate to \eqtag{gl-embedding} by the
Lie algebra automorphism of $\gog\gol(n)$ induced by
the only non-trivial automorphism of the Dynkin diagram
of $A_{n-1}$. Since there are corresponding (dual) Hopf $\ast$-algebra
automorphisms of $\Uq$ and $\Aq$, 
every statement about $\gok^\infty$-invariance can immediately 
be translated into a statement about $\gok^{-\infty}$-invariance. 
In the remainder of this paper, we shall therefore usually 
confine ourselves to proving statements for $\sigma=\infty$ and
leave it to the reader to verify the corresponding statement for
$\sigma=-\infty$.
\endth
\beginth{Theorem} Kfixed
{\thstil Let $\sigma=\pm \infty$. 
For any $\lambda \in P^+$, the subspace 
$V(\lambda)_\ks$
of $\ks$-fixed vectors in $V(\lambda)$ is at most 
one-dimensional. $V(\lambda)$ has non-zero $\ks$-fixed
vectors if and only if 
$\lambda=l(\eps_1-\eps_n)$ ($l\in \ZZ_+$).
A $\ks$-fixed vector in $V(\eps_1-\eps_n)\subset
V^\ast \ten V$ is given by $w_{\ks}= v_n^\ast \ten v_n$
($\sigma = \infty$) or $w_{\ks}= v_1^\ast \ten v_1$
($\sigma=-\infty$).
The same statements hold with $\ks$ replaced by
$\kss$.
\par}
\endth
\beginproof
It is completely straightforward to check that $w_{\gok^\infty}$
is indeed $\gok^\infty$-fixed. All the other statements concerning
 $\gok^\infty$-fixed vectors are proved in \cite{\NYM, Prop.\ 4.2}.
\endproof
\beginth{Remark} kinfty
Note that the expressions $w_{\gok^\infty}= v_n^\ast \ten v_n$ 
and $w_{\ks}= v_1^\ast \ten v_1$ are 
special cases (up to a scalar multiple) 
of \eqtag{repkfixed} ($c=1,d=0$ and $c=0, d=1$ respectively). 
In other words,
the expressions for $\ks$-fixed and $\kss$-fixed vectors in
\eqtag{repkfixed} are valid for any $c,d \geq 0$, $(c,d) \neq (0,0)$.
\endth
Specializing the notion of $\Uq$-invariant vectors
in a $\Uq$-module $W$ to the $\Uq$-action on
$\Aq$, we get the following terminology.
An element $a\in \Aq$ is {\sl left}
resp.\ {\sl right invariant} with respect to 
$u\in \Uq$ if $u\cdot a = \eps(u) a$ resp.\ 
$a\cdot u = \eps(u)a$ (cf.\ \eqtag{kfixeddef}). 

For any $-\infty \leq \sigma \leq \infty$ let us define
$$\Bq^\sigma = \Bq^{(c,d)} :=
 \{ a \in \Aq \mid a \cdot \gok^\sigma = 0 \}.
\eqno\eq{Bqdef}$$
Then $\Bq^\sigma$ is a $\ast$-subalgebra and right coideal in $\Aq$, 
as follows from Proposition \thtag{kcoideal} 
(cf.\ \cite{\DKqsf, Prop.\ 1.9}). 
It is also invariant under the left action of $\Uq$ on $\Aq$.

Recall that $W(\lambda)$ ($\lambda\in P^+$) denotes the subspace
of $\Aq$ spanned by the coefficients of the (co-)representation
$V(\lambda)$.
\beginth{Theorem} planchdecomp 
{\thstil Let $-\infty \leq \sigma \leq \infty$ be arbitrary.
The irreducible decomposition of $\Bq^\sigma$ as a
right $\Aq$-comodule resp.\ left $\Uq$-module is given by:
$$\Bq^\sigma = \bigoplus_{l\in\ZZ_+} V(l(\eps_1 - \eps_n)),
\eqno\eq{planchformula}$$
where the isotypical subspace of type $V(l(\eps_1 - \eps_n))$ 
is equal to the intersection of $\Bq^\sigma$ and $W(\lambda)$.
\par}
\endth
\beginproof
This follows from Corollary \thtag{gelfand}, Theorem 
\thtag{spherrep}, and Theorem \thtag{Kfixed}. We omit the 
details, since the proof is completely analogous to
the one given in \cite{\Noumimac, Prop.\ 4.1}.
\endproof
Our next goal will be to give a more explicit description of the
subalgebra $\Bq^\sigma \subset \Aq$.
For this we need some results from \cite{\NYM, \S 1.5} on the 
linear independence
of products of the $t_{ij}$. Let us remark that 
the commutation relations $RT_1 T_2 = T_2 T_1 R$ between the $t_{ij}$
are equivalent to the following relations:
$$\eqalignno{&t_{ij}t_{ik} = qt_{ik}t_{ij},\quad
t_{lj}t_{ij} = q^{-1} t_{ij}t_{lj},\quad
t_{lk}t_{lj}= q^{-1}t_{lj}t_{lk},\quad
t_{lk}t_{ik}= q^{-1} t_{ik}t_{lk}, & \cr
&t_{lj}t_{ik}= t_{ik}t_{lj},\quad
t_{lk}t_{ij} = t_{ij}t_{lk} - (q-q^{-1})t_{ik}t_{lj}
\quad (j<k, i<l). & \eq{tijrel}\cr}$$ 
Using this explicit form of the relations, one can construct
a basis of monomials in the $t_{ij}$ for the subalgebra of
$\Aq$ generated by the $t_{ij}$. Let $M(n,\NN)$ denote the
space of $n\times n$ matrices with non-negative integer
coefficients. For any matrix $A= (a_{ij}) \in M(n,\NN)$
define a monomial 
$t^A := t_{11}^{a_{11}} t_{12}^{a_{12}} \cdots
t_{21}^{a_{21}} t_{22}^{a_{22}} \cdots t_{nn}^{a_{nn}},$
where the factors are arranged according to 
the lexicographical order on the indices $(i,j)$.
The following theorem is proved in \cite{\NYM, Thm.\ 1.4}:
\beginth{Theorem} PBW
{\thstil The elements $t^A\in \Aq$ ($A\in M(n,\NN)$) form a
linear basis of the subalgebra $M_q\subset \Aq$ 
generated by the $t_{ij}$ ($1\leq i,j \leq n$).
\par}
\endth
\noindent
With any $A=(a_{ij})\in M(n,\NN)$ we associate a sequence
of integers $r(A)$ given by
$$r(A) := (\sum_{ij} a_{ij}, a_{11}, a_{12}, \ldots, a_{21}, a_{22},
\ldots, a_{nn})\in \NN^{n^2 +1}.$$
The lexicographic ordering on $\NN^{n^2 +1}$ then induces a total ordering
$\preceq$ on $M(n,\NN)$ in the obvious way. Note that this ordering is 
compatible with the additive structure on $M(n,\NN)$. It now
follows that any non-zero element $\phi\in M_q$ can be uniquely written 
as a linear combination of the form
$$\phi = c_A t^A + \sum_{B\prec A} c_B t^B,$$
where $c_A \neq 0$. Let us call $A =: d(\phi)$ the degree of $\phi$.
Suppose $\phi = t_{i_1j_1} \cdots t_{i_pj_p}$ is any monomial in the
$t_{ij}$. It then follows from the relations \eqtag{tijrel} and 
Theorem \thtag{PBW} that
$d(\phi) = A$, where the coefficient $a_{ij}$ 
of $A = (a_{ij})$  is defined as the number of occurrences
of $t_{ij}$ in the product $t_{i_1j_1} \cdots t_{i_pj_p}$.
It is now easy to see that  a product of non-zero elements
in $M_q$ is again non-zero and that
$$d(\phi \psi) = d(\phi) + d(\psi),\quad 0\neq \phi, \psi \in M_q.
\eqno\eq{d-add}$$
Since every element in $\Aq$ can be written as a product of
a power of $\detq^{-1}$ and an element in $M_q$ (recall that
$\detq^{-1}$ is central in $\Aq$), we get:
\beginth{Lemma} zerodiv
{\thstil The algebra $A_q$ has no (left or right) zero divisors.
\par}
\endth
Corresponding to the contragredient representation 
$V^\ast$ we have a  right $\Uq$-module $\Hom(V^\ast, \CC)$,
whose underlying vector space can be naturally identified with
$V$. The assignment $v_i \mapsto 
q^{\langle 2\rho, \eps_i\rangle} v_i^\ast$ (cf.\ \eqtag{Ssquare})
defines an intertwining operator from $\Hom(V^\ast,\CC)$
into $(V^\ast)^\circ$. The mapping 
$$V^\ast \ten \Hom(V^\ast, \CC) \to \Aq,\quad
v_i^\ast \ten v_j \mapsto t_{ji}^\ast\eqno\eq{Wdualvector}$$
is an injective $\Uq$-bimodule homomorphism. 
\beginth{Proposition} zij
{\thstil Let $-\infty \leq \sigma \leq \infty$ be arbitrary.
We put
$$x_{ij} := d t_{1i}^\ast t_{1j} + c t_{ni}^\ast t_{nj}
+ \sqrt{cd}\, t_{ni}^\ast t_{1j} + \sqrt{cd}\, t_{1i}^\ast t_{nj}
\in\Aq \quad (1\leq i,j \leq n). \eqno\eq{zdef}$$
The subspace spanned by the $x_{ij}$ ($1\leq i,j \leq n$)
is invariant under the left $\Uq$-action. The linear mapping
defined by
$$v_i^\ast \ten v_j \mapsto x_{ij}\colon  
V^\ast \ten V \to \Aq \eqno\eq{vzij}$$
is an injective operator intertwining the left $\Uq$-actions.
The $x_{ij}$ are right $\gok^\sigma$-invariant and satisfy
$x_{ij}^\ast = x_{ji}$. They generate the subalgebra
$\Bq^\sigma$. A highest weight vector for 
$V(l(\eps_1-\eps_n))\subset \Bq^\sigma$ is $x_{n1}^l$.
\par}
\endth
\beginproof
It follows from \eqtag{zdef} 
and the remarks preceding this proposition, 
 that the operator
\eqtag{vzij} intertwines the left actions of $\Uq$ on 
$V^\ast\ten V$ and $\Aq$. Hence, the subspace spanned by
the $x_{ij}$ is invariant. It is proved in
\cite{\NYM, Corol\-lary to Prop.\ 1.1} that 
$\sum_k q^{2(i-k)} t_{ik}^\ast t_{jk} = \delta_{ij}\,1\in\Aq$
($1\leq i,j\leq n$).
It follows that 
$\sum_k q^{2(n-k)}x_{kk} = c + q^{2(n-1)}d\in \Aq$.
Since $c+q^{2(n-1)}d \neq 0$, the trivial representation
occurs with non-zero multiplicity in the irreducible
decomposition of the subspace spanned by the $x_{ij}$.
On the other hand, the element $x_{n1}\in \Aq$ has 
weight $\eps_1-\eps_n$ and is annihilated by the
$L^-_{ij}$ ($i>j$), being the image of the highest weight
vector $v_n^\ast \ten v_1$.
 If we can prove that $x_{n1}$
is non-zero, it will follow from \eqtag{VVstar} that
the operator \eqtag{vzij} is injective. 
If $\sigma$ is finite, $x_{n1}$
is clearly non-zero, since $\eps(x_{n1}) =\sqrt{cd}
\neq 0$. If $\sigma=\infty$ it follows from
Lemma \thtag{zerodiv} that $x_{n1} = t^\ast_{nn} t_{n1}$
is non-zero. Next, the fact that the $x_{ij}$ are right
$\ks$-invariant follows immediately from the expression for the
$\kss$-fixed vector in \eqtag{repkfixed} and 
\eqtag{Wvector}, \eqtag{Wdualvector}. 
The property $x_{ij}^\ast = x_{ji}$
is completely trivial. It therefore remains to prove that
the $x_{ij}$ generate the whole subalgebra $\Bq^\sigma$.
Indeed, a second application of Lemma \thtag{zerodiv}
shows that $x_{n1}^l$ is non-zero. But then it is a highest weight
vector of weight $l(\eps_1-\eps_n)$. Hence, in the subalgebra
generated by the $x_{ij}$, the isotypical subspace of 
type $V(l(\eps_1-\eps_n))$ is non-zero. The assertion now
follows from Theorem \thtag{planchdecomp}. 
\endproof
\beginth{Remark} leftks
Instead of $\Bq^\sigma$, one can also consider the subalgebra
$\FSC_q^\sigma$ of left $\ks$-invariant functions. It
is a right $\Uq$-module and decomposes as: 
$$\FSC_q^\sigma = \bigoplus_{l\in\ZZ_+} V(l(\eps_1 - \eps_n))^\circ.
\eqno\eq{Cqdecomp}$$
The elements 
$y_{ij} := qd t_{i1}^\ast t_{j1} + q^{-1}c t_{in}^\ast t_{jn}
+ \sqrt{cd}\, t_{in}^\ast t_{j1} + \sqrt{cd}\, t_{i1}^\ast t_{jn}
\in\Aq$ ($1\leq i,j \leq n$)
are left $\ks$-invariant and generate the algebra $\FSC_q^\sigma$.
A highest weight vector in $V(l(\eps_1-\eps_n))^\circ$ is $y_{n1}^l$.
The assignment $v_i^\ast \ten v_j \mapsto q^{-\langle 2\rho, \eps_i
\rangle} y_{ij}$ defines an
injective intertwining operator of $(V^\ast \ten V)^\circ$
into $\FSC_q^\sigma$. 
\beginsektion 4 {A $q$-analogue of the Hopf fibration}
In this section we shall describe a more geometrically
inspired way to construct the algebra $\Bq^\sigma$ in a uniform way
for all values of $\sigma$.

Let $\Aq(\tilde S) = \Aq(\tilde S^{2n-1})$ denote the algebra 
generated by the symbols
$z_i, w_i$ ($1\leq i\leq n$) and $c,d$ subject to the
relations:
$$\eqalignno{
(i) \quad & z_i z_j = qz_jz_i,\; qw_iw_j=w_jw_i \quad 
(1\leq i<j\leq n), & \cr
(ii) \quad & w_j z_i = qz_i w_j \quad 
(1\leq i\neq j\leq n), & \eq{defrelS}\cr
(iii) \quad & w_jz_j = z_j w_j + (1-q^2) \sum_{1\leq k < j}
z_kw_k - (1-q^2)q^{-2} d \quad (1\leq j\leq n), & \cr
(iv) \quad & \sum_{k=1}^n z_kw_k = c+q^{-2}d. & \cr
}$$
As can be easily checked, these defining relations imply the following 
relations ($1\leq j\leq n$):
$$\eqalignno{
(v) \quad & cz_j = z_j c, \; cw_j = w_j c, 
\; cd =dc, & \cr
(vi) \quad & q^{-2} dz_j = z_j d, \; dw_j =q^{-2}w_jd, 
& \eq{consrelS} \cr
(vii) \quad & z_jw_j =w_jz_j + 
(1-q^{-2)} \sum_{1\leq k<j}
q^{2(j-k)} w_kz_k -(1-q^{-2})q^{2(j-1)}d, &\cr
(viii) \quad & \sum_{k=1}^n q^{2(n-k)} w_kz_k = 
c+q^{2(n-1)}d. & \cr
}$$
Note that $c$ is central, but $d$ is not.
The algebra $\Aq(\tilde S^{2n-1})$ can be 
considered as the algebra of functions on the 
{\sl total space of a 
family of quantum $(2n-1)$-spheres}. 
The case $n=2$ was considered
in \cite{\NMqthreesphere}. The algebra 
$\Aq(\tilde S^{2n-1})$ was introduced in
general in \cite{\VKproj}.

There is a unique $\ast$-structure on the algebra 
$\Aq(\tilde S^{2n-1})$ such that
$$z_j^\ast = w_j, \quad c^\ast =c, \quad d^\ast = d.$$

Let us define algebra homomorphisms
$$R\colon \Aq(\tilde S) \to 
\Aq(\tilde S)\ten \Aq, \quad
L\colon \Aq(\tilde S) \to \FSA(U(1)) \ten
\Aq(\tilde S)$$
by putting ($1\leq j\leq n$):
$$\eqalignno{R(z_j) &=\sum_{i=1}^n z_i\ten t_{ij},\;
R(w_j) = \sum_{i=1}^n w_i \ten S(t_{ji}),\;
R(c) = c\ten 1,\;R(d) = d\ten 1, & \cr
L(z_j) & = z\ten z_j,\; L(w_j) = z^{-1} \ten w_j,\;
L(c) = 1\ten c, \; L(d) = 1\ten d.& \eq{RLdef}\cr}$$
In this way, $\Aq(\tilde S)$ becomes  a two-sided 
$(\FSA(U(1)), \Aq)$-comodule  algebra.
Note that $R$ and $L$ commute with the $\ast$-operations.
Let us remark that by ``differentiating''  the right
$\Aq$-coaction one obtains a left $\Uq$-module structure
on $\Aq(\tilde S)$. 

The algebra $\Aq(\tilde S)$ can be
realized inside a suitable extension of $\Aq$ by a $q$-shift
operator.
Let $\FSC:= \CC[\alpha, \beta, \gamma, \delta]$ denote the 
free polynomial algebra in the commuting variables 
$\alpha, \beta, \gamma, \delta$. Define an algebra
automorphism $\theta\colon \FSC \to \FSC$ by 
$$\theta(\alpha)= \alpha, \quad \theta(\beta) = q\beta,
\quad \theta(\gamma) = q\gamma, \quad \theta(\delta) = \delta.
\eqno\eq{taushiftdef}$$
Let $\FSC[\theta^{\pm 1}]$ denote the subalgebra of
$\End_\CC(\FSC)$ generated by left multiplication by
$\alpha, \beta, \gamma, \delta$ and the $q$-shift operators
$\theta, \theta^{-1}$. One has the obvious algebra
isomorphism
$\FSC[\theta^{\pm 1}] \cong \CC[\theta^{\pm 1}] \ten
\CC[\alpha, \beta, \gamma, \delta],$
where the (twisted) multiplication on the tensor product
 on the right-hand side  is defined by the rule
$\theta P = \theta(P) \theta$
($P\in \CC[\alpha, \beta, \gamma, \delta]$).
There is a unique $\ast$-structure on $\FSC[\theta^{\pm 1}]$
such that
$$\theta^\ast = \theta^{-1}, \quad \alpha^\ast = \delta,\quad
\beta^\ast = -\gamma.\eqno\eq{taustardef}$$
Now consider the $\ast$-algebra $\FSC[\theta^{\pm 1}] \ten
\Aq$, the multiplication and the $\ast$-structure on the tensor
product being defined as usual. We define a two-sided
$(\FSA(U(1)), \Aq)$-comodule structure on 
$\FSC[\theta^{\pm 1}] \ten \Aq$ with structure mappings
$$R\colon \FSC[\theta^{\pm 1}] \ten \Aq \to
\FSC[\theta^{\pm 1}] \ten \Aq \ten \Aq, \quad 
R:= \id \ten \Delta,$$
$$L\colon \FSC[\theta^{\pm 1}] \ten \Aq \to
 \FSA(U(1))\ten \FSC[\theta^{\pm 1}] \ten \Aq,\quad
L(\theta) := z^{-1} \ten \theta,$$ 
$L$ acting trivially on all the other generators.
\beginth{Theorem} Phi-embedding
{\thstil There is a unique algebra homomorphism
$\Phi\colon \Aq(\tilde S) \to\FSC[\theta^{\pm 1}] \ten \Aq$
such that $(1\leq j\leq n$):
$$\Phi(z_j) = \theta^{-1} (\gamma t_{1j} + \delta t_{nj}), \;
\Phi(w_j) = (-\beta t_{1j}^\ast + \alpha t_{nj}^\ast)\theta,\;
\Phi(c) = \alpha\delta,\; \Phi(d) = -\beta\gamma. 
\eqno\eq{Phidef}$$
The mapping $\Phi$ is injective, commutes with the $\ast$-structures,
and intertwines the two-sided $(\FSA(U(1)), \Aq)$-comodule
structures.
\par}
\endth
\beginproof
All assertions are obvious except the existence and the injectivity.
To check the existence
we need some information about the commutation relations 
between the $t_{ij}$ and the $t_{ij}^\ast$ in $\Aq$. Put 
$\check R = PR$ with $P\colon V\ten V\to V\ten V$ the permutation
operator. The commutation relations between the $t_{ij}$
 can be rewritten as
$\check R T_1 T_2 = T_1 T_2 \check R$. Using $T S(T) = S(T) T = I$
one deduces the identity
$S(T)_1 \check R T_1 = T_2 \check R S(T)_2.$
Making this explicit one gets
($1\leq i,j,k,l \leq n$):
$$q^{\delta_{jl}} t_{li}^\ast t_{jk} + (q-q^{-1}) \delta_{jl}
\sum_{m<j} t_{mi}^\ast t_{mk} =
q^{\delta_{ik}} t_{jk} t_{li}^\ast + (q-q^{-1}) \delta_{ik}
\sum_{m>i} t_{jm} t_{lm}^\ast. \eqno\eq{ttstar}$$
Using these relations it is straightforward to check that $\Phi$ is
well-defined. To prove that $\Phi$ is injective, first observe 
that the algebra $\Aq(\tilde S)$ has a natural grading preserved by
the $\Uq$-action 
such that $\deg(z_i) = \deg(w_i) = 1$ ($1\leq i\leq n$) and
$\deg(c) = \deg(d) = 2$. Suppose the $\Phi(w_n^m z_1^l)$ are non-zero. 
Then they are highest weight vectors of weight $l\eps_1 - m\eps_n$
with respect to the $\Uq$-action, and so are the 
$w_n^m z_1^l\in\Aq(\tilde S)$.
Let us denote the $\Uq$-submodule
generated by $w_n^m z_1^l$ ($l,m$ fixed) by $V(l\eps_1 -m\eps_n)$.
The elements $w_n^m z_1^l c^r d^s$ ($l,m$ fixed, $r,s\geq 0$ arbitrary)
are linearly independent in $\Aq(\tilde S)$, since their images under
$\Phi$ are in $\FSC[\theta^{\pm 1}]\ten \Aq$.
This shows that multiplication defines an injective linear mapping
$$\bigoplus_{l,m\geq 0} V(l\eps_1 -m\eps_n) \ten \CC[c,d]
\hookrightarrow \Aq(\tilde S)\eqno\eq{injec}$$ intertwining the
$\Uq$-actions. On the other hand, by \eqtag{defrelS}
multiplication defines a surjective graded linear mapping
$$\CC_q[w_1, \ldots, w_n] \ten \CC_q[z_1, \ldots, z_n]\ten
\CC[d]\twoheadrightarrow \Aq(\tilde S)\eqno\eq{surj}$$
intertwining the $\Uq$-actions. Here the algebras on the 
left-hand side are the well-known polynomial algebras in $q$-commuting
variables (cf.\ \eqtag{defrelS}) with the natural $\Uq$-symmetry and
grading. Their respective  irreducible decompositions are given by
$\CC_q[z_1, \ldots, z_n] = \bigoplus_{l\geq 0} V(l\eps_1)$
and $\CC_q[w_1, \ldots, w_n] = \bigoplus_{m\geq 0} V(-m\eps_n)$.
For any $N\geq 0$, the mappings
\eqtag{injec} and \eqtag{surj} place lower and upper bounds
respectively on the dimension of the homogeneous component
in $\Aq(\tilde S)$ of degree $N$. An easy counting argument
using the well-known fact that the tensor product 
$V(l\eps_1) \ten V(-m\eps_n)$ decomposes as the direct sum
$\bigoplus_{i=0}^{l\wedge m}V((l-i)\eps_1 - (m-i)\eps_n)$
shows that these upper and lower bounds are actually
equal. This implies that the mapping \eqtag{injec} is surjective.
Obviously, the restriction of $\Phi$ to
each of the components $V(l\eps_1 -m\eps_n) \ten \CC[c,d]$ 
in \eqtag{injec} is injective,
the $\Uq$-modules $V(l\eps_1 -m\eps_n)$ being irreducible.
The injectivity of $\Phi$ therefore follows from the
following lemma.
\endproof
\beginth{Lemma} wznonzero
{\thstil The elements 
$$\Phi(w_n^m z_1^l) = 
\left \{ (-\beta t_{1n}^\ast + \alpha t_{nn}^\ast)\theta \right \}^m
\left \{ \theta^{-1} (\gamma t_{11} + \delta t_{n1}) \right \}^l \quad
(l,m \geq 0)$$
in the algebra $\FSC[\theta^{\pm 1}]\ten \Aq$ are non-zero.
\par}
\endth
\beginproof
We are going to apply Theorem \thtag{PBW}, but for this
we need a more explicit description of 
the $\ast$-operation. For any $1\leq i\leq n$ 
denote by $(i_1, \cdots, i_{n-1})$
the sequence of length $n-1$ obtained from $(1, \ldots, n)$ by
deleting $i$. Define the quantum principal
minor $\xi^i_j\in \Aq$ by
$$\xi^i_j := \sum_{w\in \goS_{n-1}} (-q)^{l(w)}t_{i_{w(1)}j_1} \cdots
t_{i_{w(n-1)}j_{n-1}}\quad
(1\leq i,j \leq n).\eqno\eq{minordef}$$
Here $\goS_{n-1}$ is the permutation group on $n-1$ letters,
and $l(w)$ denotes the length of $w$.
It is well-known (cf.\ \cite{\RTF, Thm.\ 4}, \cite{\NYM, (1.24)})
that the $\ast$-operation is given on the
generators $t_{ij}$ ($1\leq i,j \leq n$) by
$$t_{ij}^\ast = (-q)^{j-i} \detq^{-1}\xi^i_j
\quad (1\leq i,j \leq n). \eqno\eq{starminor}$$
Recall the notion of degree for a monomial in the $t_{ij}$ (cf.\ section
3). On the one hand, $d(t_{n1}) \prec d(t_{11})$. On the other hand,
it is easy to see from \eqtag{minordef} that $d(\xi^1_n) \prec d(\xi^n_n)$. 
Let us put $A:= d(t_{11})$ and $B:= d(\xi^n_n)$.
Since $t_{11} t_{n1} = q t_{n1} t_{11}$ and 
$t_{nn}^\ast t_{1n}^\ast = q t_{1n}^\ast t_{nn}^\ast$, we can
expand $\Phi(w_n^m z_1^l)$ as 
$$
\sum_{M = (m_1, m_2, l_1, l_2)} c_M P_M(\theta, \alpha, \beta, \gamma, 
\delta) \ten \detq^{-m_1-m_2} (\xi^1_n)^{m_1} (\xi^n_n)^{m_2}
t_{11}^{l_1} t_{n1}^{l_2},$$
where the sum runs over all quadruples $M = (m_1, m_2, l_1, l_2)$ of
non-negative integers such that $m_1 + m_2 = m$ and $l_1 + l_2 = l$,
the $c_M$ are non-zero constants, and the 
$P_M(\theta, \alpha, \beta, \gamma, \delta)$ are certain (commutative) 
monomials in the variables $\theta, \alpha, \beta, \gamma, \delta$.
A moment's consideration shows that for all $(m_1, m_2, l_1, l_2)$
$$d((\xi^1_n)^{m_1} (\xi^n_n)^{m_2} t_{11}^{l_1} t_{n1}^{l_2}) \preceq
mB + lA,$$
and that equality holds if and only if $(m_1, m_2) = (0, m)$ and
$(l_1, l_2) = (l, 0)$. The corresponding coefficient $c_{(0,m,l,0)}$ equals
$q^k \alpha^m \gamma^l$ for a certain $k\in \ZZ$. 
This proves that $w_n^m z_1^l$ is non-zero.
\endproof\noindent
As a corollary of the proof of Theorem \thtag{Phi-embedding}
we have:
\beginth{Corollary} ASdecomp
{\thstil The isotypical decomposition of $\Aq(\tilde S)$
with respect to the $\Aq$-coaction is:
$$\Aq(\tilde S) \cong \bigoplus_{l,m\geq 0} 
V(l\eps_1 -m\eps_n) \ten \CC[c,d].$$
Here $V(l\eps_1 -m\eps_n)$ is the $\Aq$-subcomodule of $\Aq(\tilde S)$
generated by the highest weight vector $w_n^m z_1^l$. The
isomorphism $\cong$ is given (from right to left) by multiplication.
\par}
\endth\noindent
Let us now consider $U(1)$-invariance in $\Aq(\tilde S)$.
Recall that in the classical case dividing out the sphere
$S^{2n-1}\subset \CC^n$ by its natural $U(1)$-action
yields the complex projective space $\CC\PP^{n-1}$ (Hopf
fibration).
We define 
$$\Aq(\widetilde{\CC\PP}{}^{n-1}):=
\{a\in \Aq(\tilde S^{2n-1}) \mid L(a) = 1\ten a \}.
\eqno\eq{defprojcd}$$ 
Clearly, $\Aq(\widetilde{\CC\PP}{}^{n-1})$ is a $\ast$-subalgebra
and $\Aq$-subcomodule of $\Aq(\tilde S)$. 
It follows immediately from Corollary \thtag{ASdecomp} that
$\Aq(\widetilde{\CC\PP}{}^{n-1})$ is generated by the elements
$\tilde x_{ij} := w_i z_j \in \Aq(\tilde S)$ and $c,d$. 
Note that both $c$ and $d$ are central in
$\Aq(\widetilde{\CC\PP}{}^{n-1})$. 
Let $\Aq(\CC\PP_q^{n-1}(c,d))$ denote the algebra obtained 
from $\Aq(\widetilde{\CC\PP}{}^{n-1})$ by specialization of 
$c,d$ to real non-negative values not both equal to zero.
It has a natural right $\Aq$-comodule
structure with an irreducible  decomposition given by
$$\Aq(\CC\PP_q^{n-1}(c,d)) = \bigoplus_{l\geq 0} 
V(l(\eps_1 -\eps_n)).\eqno\eq{CPdecomp}$$
The algebras $\Aq(\widetilde{\CC\PP}{}^{n-1})$ and its specialization
already appeared in \cite{\VKproj}. For  $n=2$, the algebras
$\Aq(\CC\PP^{n-1}(c,d))$ ($c,d\geq 0$)  coincide with
the algebras of functions on the quantum spheres
introduced in \cite{\Pod}.
\beginth{Theorem} special
{\thstil Let $c_0,d_0 \geq 0$ be real numbers not both equal to zero.
There is a unique algebra homomorphism
$$\Aq(\widetilde{\CC\PP}{}^{n-1}) \to \Bq^{(c_0,d_0)},
\quad \tilde x_{ij} \mapsto x_{ij}, \; c \mapsto c_0,
\; d\mapsto d_0.\eqno\eq{specialmap}$$
This mapping induces a $\ast$-algebra isomorphism
of $\Aq(\CC\PP_q^{n-1}(c_0,d_0))$ onto $\Bq^{(c_0,d_0)}$
intertwining the $\Aq$-coactions.
\par}
\endth
\beginproof
Define an algebra homomorphism $\phi\colon 
\CC[\alpha, \beta, \gamma, \delta] \to \CC$ by the rule
$$\phi(\alpha) = \phi(\delta) = \sqrt{c_0},
\quad \phi(\beta) = -\sqrt{d_0}, \quad \phi(\gamma) = \sqrt{d_0}.$$
Obviously, $\phi$ maps $c,d\in \CC[\alpha,
\beta, \gamma, \delta]$ onto $c_0,d_0 \in \CC$.
Since $\Phi(\Aq(\widetilde{\CC\PP}{}^{n-1}))$ is contained in
the subspace $\CC[\alpha, \beta, \gamma, \delta] \ten
\Aq$ of $\FSC[\theta^{\pm 1}] \ten \Aq$, the composition
$(\phi\ten \id) \circ \Phi$ is well-defined on
 $\Aq(\widetilde{\CC\PP}^{n-1})$.  It is a 
$\ast$-algebra homomorphism into $\Aq$ mapping 
$\tilde x_{ij} = w_iz_j$ onto
$$\eqalign{(-\phi(\beta) t_{1i}^\ast &+ \phi(\alpha) t_{ni}^\ast)
(\phi(\gamma) t_{1j} + \phi(\delta) t_{nj}) \cr
& = d_0 t_{1i}^\ast t_{1j} + c_0 t_{ni}^\ast t_{nj}
+ \sqrt{c_0 d_0}\, t_{ni}^\ast t_{1j} + 
\sqrt{c_0 d_0}\, t_{1i}^\ast t_{nj}.\cr}$$
This proves the existence of \eqtag{specialmap}. The remaining
assertions now follow trivially from \eqtag{CPdecomp}, Theorem
\thtag{planchdecomp}, Prop.\ \thtag{zij} and the definitions.
\endproof
We shall call either $\Aq(\CC\PP^{n-1}(c,d))$ or $\Bq^{(c,d)}$ 
the algebra of representative functions on the
{\sl quantum projective space} $\CC\PP_q^{n-1}(c,d)$.
\beginth{Corollary} tilde-eps
{\thstil There exists a unique $\ast$-algebra homomorphism
$$\tilde\eps\colon \Aq(\CC\PP^{n-1}(c,d)) \to \CC$$
such that 
$\tilde\eps(x_{11}) = d$, $\tilde\eps(x_{nn}) = c$,
$\tilde\eps(x_{1n}) = \sqrt{cd}$,
$\tilde\eps(x_{n1}) = \sqrt{cd}$, and 
$\tilde\eps(x_{ij})=0$ in all other cases.
\par}
\endth
\beginproof
Simply take $\tilde\eps$ equal to the composition of the
isomorphism described in Theorem \thtag{special} and the 
counit $\eps\colon \Aq \to \CC$.
\endproof\noindent
The $\ast$-homomorphism $\tilde\eps$ can be regarded as a 
``classical'' point in the quantum projective space
$\CC\PP_q^{n-1}(c,d)$ (cf.\ \cite{\DKqsf, Prop.\ 1.1}).
\beginsektion 5 {Zonal spherical functions}
Let us fix parameters $-\infty \leq\si,\tau\leq \infty$. We define
the $\ast$-subalgebra of {\sl $(\si,\tau)$-biinvariant
functions} in $\Aq$ as 
$$\Hst := \{ a\in \Aq \mid
\gok^\si \cdot a = 0 \quad \hbox{and} 
\quad a\cdot \gok^\tau = 0\},\eqno\eq{Hdef}$$
An element $a\in \Aq$ is called {\sl 
$(\si,\tau)$-spherical} if it is $(\si,\tau)$-biinvariant
and contained in $W(\lambda)$ for some $\lambda\in P^+$.
\beginth{Proposition} Hdecomp
{\thstil If we put $\Hst(\lambda) := \Hst\cap W(\lambda)$, then 
$$\Hst = \bigoplus_{l\in\ZZ_+} \Hst(l(\eps_1 -\eps_n)),
\eqno\eq{Hdecomp}$$
and each of the spaces $\Hst(l(\eps_1 -\eps_n))$ is
one-dimensional. 
\par}
\endth
\beginproof
This is a direct consequence of Theorems \thtag{planchdecomp},
\thtag{Kfixed}, \thtag{spherrep}, and Corollary \thtag{gelfand}.
\endproof
Let us now suppose that $\si,\tau$ are finite.
It follows immediately from \eqtag{repkfixed} and 
Proposition \thtag{zij} that the direct sum $\Hst(0) \oplus \Hst(\eps_1
-\eps_n)$ is spanned by the unit element $1\in \Aq$
and the element 
$$\xst := {1\over 2} (\hat x_{1n} + \hat x_{n1} + q^{\si +1}
\hat x_{11} + q^{-\si-1} \hat x_{nn} -
(q^{\si+\tau+1} + q^{-\si-\tau-1}))
\in \Aq. \eqno\eq{kbiinvar}$$
Here we take ($q^\tau = \sqrt{d/c}$):
$$\hat x_{ij} := q^\tau t_{1i}^\ast t_{1j} + 
q^{-\tau} t_{ni}^\ast t_{nj}
+ t_{ni}^\ast t_{1j} + t_{1i}^\ast t_{nj}
\in\Aq \quad (1\leq i,j \leq n). \eqno\eq{defztau}$$
Note that the $\hat x_{ij}$ differ from the $x_{ij}$
(cf.\ \eqtag{zdef}) by a scalar multiple.
One has $(\xst)^\ast = \xst$. The following lemma is obvious:
\beginth{Lemma} zrestriction
{\thstil Let $\si,\tau$ be finite. 
Under the restriction mapping 
${}_{|\TT}\colon \Aq \longrightarrow \AT$ we have:
$$\hat x_{11} \mapsto q^\tau, \quad \hat x_{nn} 
\mapsto q^{-\tau}, \quad
\hat x_{1n} \mapsto z_1^{-1} z_n, \quad \hat x_{n1} 
\mapsto z_1z_n^{-1},$$
and all the other $x_{ij}$ are mapped onto $0$. 
In particular, the image of $\xst$ equals 
${1\over 2}(z_1z_n^{-1} + z_1^{-1}z_n)$.
\par}
\endth
\beginth{Proposition} imageH
{\thstil Let $\si,\tau$ be finite. The algebra
$\Hst$ is generated by $\xst$ and hence commutative.
If we put $z := z_1z_n^{-1}\in \AT$ then
the restriction of the mapping ${}_{|\TT}$ to 
$\Hst\subset \Aq$
is an injective $\ast$-algebra homomorphism 
onto the polynomial algebra 
$\Hst_{|\TT}:=\CC[{1\over 2} (z+z^{-1})] 
\subset \AT$. 
\par}
\endth
\beginproof
Recall that $\hat x_{n1}^l\in \Aq$ ($l\in \ZZ_+$) 
is a (non-zero) highest weight vector of
weight $l(\eps_1-\eps_n)$ (cf.\ Proposition \thtag{zij}). 
Hence the component of $(\xst)^l$ in
$\Hst(l(\eps_1-\eps_n))$ is non-zero. This implies
that $\xst$ generates the algebra $\Hst$.
By Lemma \thtag{zrestriction}
the restriction of $(\xst)^l$ 
to $\TT$ is equal to a non-zero scalar multiple of
$z_1^lz_n^{-l}$ plus some lower order terms w.r.t.\ 
the lexicographic ordering on the monomials in the 
$z_i^{\pm 1}$. Therefore the images of
the $(\xst)^l$ ($l\in \ZZ_+$) in $\AT$ are linearly 
independent, which completes the proof.
\endproof 
For every $l\in\ZZ_+$, any non-zero element of 
$\Hst(l(\eps_1-\eps_n))$ can be expressed as a polynomial of 
degree $l$ in $\xst$.  In order to identify 
these polynomials, we study the action of
the following {\sl Casimir operator} (cf.\ 
\cite{\RTF}, \cite{\Noumimac}) on $\Hst$:
$$C := \sum_{ij} q^{2(n-i)}\Lp_{ij}S(\Lm_{ji})\in 
\Uq.\eqno\eq{Cdef}$$
The element $C\in\Uq$ is central and it acts as a scalar
on each subspace $W(\lambda)\in\Aq$ ($\lambda\in P^+$).
The corresponding eigenvalue is given by
$$\chi_\lambda(C) :=\sum_{k=1}^n q^{2(\lambda_k + n-k)}.
\eqno\eq{CscalarW}$$
Since $C$ is central, the left action of
$C$ on $\Aq$ preserves the subalgebra $\Hst$. Hence
it acts as a scalar on each subspace 
$\Hst(l(\eps_1-\eps_n))$ $(l\in \ZZ_+$):
 $$\chi_l(C) := \chi_\lambda(C)= q^{2(l+n-1)} + q^{-2l}
+  {q^2 -q^{2n-2}\over 1-q^2},
\quad \lambda = l(\eps_1 -\eps_n).
\eqno\eq{CscalarH}$$
It follows from Proposition \thtag{imageH} that 
there is a uniquely determined linear operator
$D\colon \Hst_{|\TT} \to \Hst_{|\TT}$ 
(called {\sl the
radial part} of the Casimir operator $C$) such 
that on $\Hst$ we have 
$${}_{|\TT} \circ C = D\circ {}_{|\TT},
\eqno\eq{radcommut}$$
where the symbol $C$ denotes the 
left action of the element
$C\in\Uq$ on $\Hst\subset \Aq$.
Let us define a linear operator $T_{q,z} \colon
\CC[z^{\pm 1}] \to \CC[z^{\pm 1}]$ (called $q$-shift operator) 
by 
putting $T_{q,z} f(z) := f(qz)$. Recall the notation
$A(z;q)$ from \eqtag{Azqdef}.
\beginth{Theorem} radialpart
{\thstil Let $\si,\tau$ be finite. 
The radial part $D\colon \Hst_{|\TT} \to \Hst_{|\TT}$ 
of the Casimir operator $C$ is
equal to the following second-order $q$-difference
operator:
$$D_{AW} := A(z;q^2)(T_{q^2,z} - \id)  + 
A(z^{-1};q^2)(T_{q^{-2},z} -\id)
+ {1-q^{2n} \over 1-q^2}\cdot\id,
\eqno\eq{radD}$$
with parameters $a,b,c,d$ given by
$$a = -q^{\si + \tau + 1},
\quad b= -q^{-\si - \tau + 1},
\quad c = q^{\si - \tau + 1}, \quad
d = q^{-\si + \tau + 2(n-2) + 1}.$$
\par}
\endth
\beginproof
The proof of this theorem is rather long and will
be deferred to section 6.
The  case $n=2$ was essentially 
proved in \cite{\Kaw, Lemma 5.1}.
\endproof
\beginth{Theorem} spherAW
{\thstil Let $\si,\tau$ be finite. The 
$(\si,\tau)$-spherical functions in  
$\Hst(l(\eps_1-\eps_n))$ $(l\in\ZZ_+)$ 
are spanned by
$$p_l(\xst; -q^{\si + \tau + 1}, 
-q^{-\si - \tau + 1},
q^{\si - \tau + 1}, 
q^{-\si + \tau + 2(n-2) + 1} \mid q^2),\eqno\eq{parameter}$$
where $p_l$ is an Askey-Wilson polynomial.
\par}
\endth
\beginproof
Recall that Askey-Wilson polynomials are characterized
as the polynomial solutions of the second-order
$q$-difference equation \eqtag{qdiffeq}. Let $\phi^l$ be
a non-zero element in $\Hst(l(\eps_1-\eps_n))$. Then
$D\,\phi^l_{|\TT} = \chi_l(C)\,\phi^l_{|\TT}$
  by definition of $D$.
 By comparing \eqtag{radD} and 
\eqtag{CscalarH} with \eqtag{qdiffeq} one sees that 
$\phi^l_{|\TT}$ is a scalar multiple of the 
 Askey-Wilson polynomial of degree $l$ in the variable
${1\over 2} (z+z^{-1})$. The theorem now follows by
Proposition \thtag{imageH}.
\endproof
For $n=2$ this result agrees 
with \cite{\Kaw, Theorem 5.2}.
As a corollary of Theorem \thtag{spherAW} we have:
\beginth {Proposition} haar
{\thstil Let $\sigma, \tau$ be finite.
Let $d m_q(z)= d m_{a,b,c,d;q}(z)$ denote the measure
defined in the right-hand side of \eqtag{orthaw-rel} and normalized in
such a way that $\int d m_q(z) = 1$ (cf.\ \eqtag{awnorm}). Then
the Haar functional $h\colon \Aq \to \CC$
is given on $\Hst$ by
$$h(P(\xst)) = \int P\left ({z+ z^{-1} \over 2}\right ) 
d_{q^2} m(z) \quad (P\in \CC[x]),$$
where the parameters $a,b,c,d$ are defined as in 
Theorem \thtag{radialpart}.
\par}
\endth
Let us now suppose $\sigma$ finite, $\tau=\pm\infty$. 
As follows from remark \thtag{leftks} and \eqtag{repkfixed}, 
the direct
sum $\FSH^{(\sigma, \pm\infty)}(0)
\oplus \FSH^{(\sigma, \pm\infty)}((\eps_1
-\eps_n))\subset \Aq$ is spanned by the 
unit element $1\in \Aq$
and
$$x^{(\sigma, \infty)} := q^{\sigma-1}\hat y_{nn} -q^{-2}
\quad \hbox{resp.}\quad
x^{(\sigma, -\infty)} := q^{-\sigma-1}\hat y_{11} -1, 
\eqno\eq{defsigmainf}$$
where $\hat y_{11},\hat y_{nn}\in\FSC_q^\sigma$ 
are given by (cf.\ remark \thtag{leftks}):
$$\hat y_{ij} := t_{in}^\ast t_{j1} + t_{i1}^\ast t_{jn} +
q^{\sigma + 1} t_{i1}^\ast t_{j1} + q^{-1-\sigma} t_{in}^\ast
t_{jn}. \eqno\eq{tildeyij}$$
Note that $\hat y_{ij}$ and $y_{ij}$ ($q^\sigma = \sqrt{{d\over c}}$)
differ by a scalar multiple.
\beginth{Proposition} sigmainfHpol
{\thstil Let $\sigma$ be finite. The algebra 
$\FSH^{(\sigma,\pm \infty)}$ is polynomial in 
$x^{(\sigma, \pm\infty)}$ and hence
commutative. 
\par}
\endth
\beginproof
Suppose $\tau=\infty$.
We may view $x^{(\sigma, \infty)}$ as a left $\ks$-invariant
function in the algebra $\Bq^\infty$.
Recall (cf.\ Proposition \thtag{zij}) that $(t_{nn}^\ast
t_{n1})^l\in \Bq^\infty$ is a non-zero highest weight vector of weight
$l(\eps_1-\eps_n)$ under the
left action of $\Uq$. Since $(t_{nn}^\ast t_{n1})^l$ is (up to
a non-zero scalar multiple) the highest weight component of
$(x^{(\sigma, \infty)})^l$, the component of
$(x^{(\sigma, \infty)})^l$ in 
$\FSH^{(\sigma, \infty)}(l(\eps_1-\eps_n))$ is non-zero, 
whence the statement for $\tau=\infty$. The case $\tau=-\infty$ is
handled in a similar way.
\endproof
By applying the limit transition from Askey-Wilson
polynomials to big $q$-Jacobi polynomials 
(cf.\ Proposition \thtag{limitbig}) 
we obtain:
\beginth{Theorem} spherbigqJ
{\thstil Let us suppose that $\sigma$ is finite.
The $(\sigma,\pm\infty)$-spherical functions in
$\FSH^{(\sigma,\pm\infty)}(l(\eps_1-\eps_n))$ $(l\in \ZZ_+)$
are spanned by
$$P_l^{(n-2,0)}(x^{(\sigma,\infty)}; 
q^{2\sigma},1\colon q^2) \quad
\hbox{resp.} \quad
P_l^{(n-2,0)}(x^{(\sigma, -\infty)}; 
q^{-2\sigma +2(n-2)},1\colon q^2),$$
where $P_l^{(n-2,0)}$ is a big $q$-Jacobi polynomial.
\par}
\endth
\beginproof
Let us for instance consider the case $\tau\to\infty$. Take $a:= q^\tau$,
$c:= q^{2\sigma}$, $d:= 1$, $\alpha := n-2$, $\beta:= 0$. 
Then $2a(c/d)^{1\over 2}x^{(\sigma,\tau)}/q \to x^{(\sigma, \infty)}$
when $\tau\to\infty$. Moreover, the parameter values in 
Proposition \thtag{limitbig} agree with those in Theorem \thtag{spherAW},
taking into account the fact 
that the Askey-Wilson polynomials $r_l(x; a,b,c,d \mid q)$
are symmetric in the parameters $a,b,c,d$ up to a
scalar multiple.
\endproof\noindent
For $n=2$ the result Theorem \thtag{spherbigqJ} 
agrees with \cite{\NMbigJ, Theorem 3}
and \cite{\Kaw, Theorem 6.2}.

Let now $\si=\pm\infty$, $\tau=\pm\sigma$. 
Then the direct
sum $\FSH^{(\si,\pm\sigma)}(0)\oplus 
\FSH^{(\si,\pm\sigma)}((\eps_1 -\eps_n))\subset 
\Aq$ is spanned by the unit element $1\in \Aq$
and $$\eqalignno{&x^{(\infty,\infty)} := 
q^{-2} t_{nn}^\ast t_{nn}-q^{-2}, \quad
x^{(-\infty,-\infty)} := t_{11}^\ast t_{11} -1,& \cr 
&x^{(\infty,-\infty)} := 
q^{-2(n-1)} t^\ast_{1n}t_{1n}, \quad
x^{(-\infty,\infty)} := 
t^\ast_{n1} t_{n1} & \eq{defxinfinf}\cr}$$
respectively.
\beginth {Proposition} infinfHpol
{\thstil Let $\sigma=\pm\infty$. The algebra 
$\FSH^{(\si,\pm\sigma)}$ is polynomial in 
$x^{(\si,\pm\sigma)}$.
\par}
\endth
\beginproof
This can be proved using Theorem \thtag{PBW} and
the description of the $\ast$-operation in terms of 
quantum principal minors (cf.\ \eqtag{minordef}, \eqtag{starminor}).
In fact, an easy argument by comparing degrees of 
monomials shows that the $(x^{(\si,\pm\sigma)})^m$ ($1\leq m\leq l$)
are linearly independent. On the other hand, they are clearly
contained in the subspace $\bigoplus_{1\leq m \leq l} 
\FSH^{(\si,\pm\si)}(m(\eps_1-\eps_n))$, which has dimension
$l$ by Proposition \thtag{Hdecomp}. This proves the assertion.
\endproof 
\beginth{Theorem} spherlittleqJ
{\thstil Let $\si=\pm\infty$. 
The $(\si,\pm\si)$-spherical functions in
$\FSH^{(\si,\pm\si)}(l(\eps_1-\eps_n))$ $(l\in \ZZ_+)$
are spanned by
$$p_l^{(n-2,0)}(x^{(\si,\si)}\colon q^2)
\quad\hbox{resp.}\quad 
p_l^{(0,n-2)}(x^{(\si,-\si)}\colon q^2),$$
where $p_l^{(n-2,0)}$ and
$p_l^{(0,n-2)}$ are little $q$-Jacobi polynomials.
\par}
\endth
\beginproof
We apply Proposition \thtag{limitlittle} to Theorem \thtag{spherAW}.
Suppose for instance $\sigma = \tau$. In the case $\sigma \to \infty$
we take $a:= q^\sigma$, $\alpha := 0$, $\beta := n-2$.
In the case  $\sigma \to -\infty$ we take $a:= q^{-\sigma}$,
$\alpha := 0$, $\beta := n-2$. 
As is easily checked, in both  cases 
$2a^2 x^{(\sigma,\tau)}/q \to x^{(\pm \infty, \pm \infty)}$
when $\sigma \to \pm \infty$.
Also, the parameter values in Proposition \thtag{limitlittle}
agree with those in Theorem \thtag{spherAW}. The result now follows. 
\endproof\noindent 
The result for the case $(\infty, \infty)$ is 
in agreement with \cite{\NYM, Theorem 4.7}. 
For $n=2$ see also \cite{\Kaw, Theorem 6.4}.
\beginth{Remark} jacrad
It is easy to deduce from Theorem \thtag{spherbigqJ} 
that the restriction of the 
Casimir operator to the subalgebra
$\FSH^{(\sigma, \pm \infty)}$ ($\sigma$ finite) 
is essentially
equal to the second-order $q$-difference operator diagonalized
by the big $q$-Jacobi polynomials (cf.\ section 2). A similar
statement holds for little $q$-Jacobi polynomials.
\beginsektion 6 {Computation of the radial part}
Throughout this section we assume that $\sigma, \tau$ are finite. We
start by establishing some notation and recalling a few facts.

From the fact  that multiplication $\Aq \ten \Aq \to \Aq$
is a $\Uq$-bimodule homomorphism one easily deduces that
$$\langle u, v\cdot a \cdot w \rangle = \langle wuv , a\rangle
\quad (u,v,w\in\Uq,\, a\in\Aq).\eqno\eq{dualaction}$$
Hence,  an element $a\in\Aq$ is 
$(\sigma, \tau)$-biinvariant if and only if 
$\langle \Uq\gok^\sigma + \gok^\tau\Uq, a \rangle = 0$.
For $u,v\in\Uq$ we write $u\sim v$ if $u-v \in \Uq \gok^\sigma +
\gok^\tau \Uq$.

Recall that the duality between $\Uq$ and $\Aq$ allows us to identify
$\Aq$ with a vector subspace of the linear dual $\Uq^\ast$. 
Restricting to $\Hst$ one has the following natural commutative diagram
$$\CD 
\Hst @>>> \FSA(\TT) \\
@VVV      @VVV \\
(\Uq/\Uq\gok^\sigma + \gok^\tau\Uq)^\ast @>>> \Uh^\ast.
\endCD \eqno\eq{comdiam}$$
Here the vertical mappings and the upper horizontal mapping are
injective.
The goal of this section is to prove the following theorem:
\beginth{Theorem} qhC
{\thstil Let $C$ denote the Casimir operator defined in
\eqtag{Cdef}. If we put $\lambda := \langle h, \eps_1 -\eps_n \rangle$
$(h\in P^\ast)$, then for generic $h\in P^\ast$:
$$q^hC \sim A(q^\lambda;q^2)(q^{h+2\eps_1} - q^h) +
A(q^{-\lambda};q^2)(q^{h-2\eps_1}-q^h) + {1-q^{2n}\over 1-q^2}\cdot q^h
\eqno\eq{qhCformula}$$
with $A(z;q)$ as defined in \eqtag{Azqdef} and $a,b,c,d$ given by
$$a = -q^{\si + \tau + 1},
\quad b= -q^{-\si - \tau + 1},
\quad c = q^{\si - \tau + 1}, \quad
d = q^{-\si + \tau + 2(n-2) + 1}.$$
Here ``generic'' means that the denominators of 
$A(q^\lambda;q^2)$ and $A(q^{-\lambda};q^2)$ are non-zero.
\par}
\endth
\beginproof
Let us first show that Theorem \thtag{qhC} implies 
Theorem \thtag{radialpart}.
Write $x:= \xst\in\Hst$ and recall 
(cf.\ \thtag{zrestriction}) that 
$x_{|\TT} = {1\over 2}(z+ z^{-1})$. Let $P_l$ be the unique
monic polynomial (of degree $l$) such that $P_l(x) \in 
\Hst(l(\eps_1-\eps_n))$. Write $R_l(z):= P_l({1\over 2}(z+ z^{-1}))$
for the corresponding Laurent polynomial in $z$.
Denote  the $q$-difference operator defined
in the right-hand side of \eqtag{radD} by $D_{AW}$.
By definition of the radial part (cf.\ \eqtag{radcommut})
 and by \eqtag{comdiam} it suffices to prove that for generic
$h\in P^\ast$:
$$\langle q^h, C\cdot P_l(x)\rangle = 
\langle q^h, D_{AW}\cdot R_l(z)\rangle
\quad (h\in P^\ast). \eqno\eq{CD}$$
On the one hand  we have
$$\eqalignno{D_{AW}\cdot R_l(z) &= 
A(z;q^2)(R_l(q^2z) - R_l(z))& \cr
&-A(z^{-1};q^2) (R_l(q^{-2}z)- R_l(z))
+ {1-q^{2n}\over 1-q^2} R_l(z).& \eq{hD}\cr}$$
On the other hand, using  \eqtag{dualaction} and
\eqtag{qhCformula}
one gets:
$$\eqalignno{\langle q^h, C\cdot P_l(x)\rangle &= 
\langle q^h C, P_l(x)\rangle = \langle A(q^\lambda;q^2)(q^{h + 2\eps_1}
- q^h) &\eq{hC} \cr
&\quad -A(q^{-\lambda};q^2) (q^{h  -2\eps_1}- q^h)
+ {1-q^{2n}\over 1-q^2} q^h, R_l(z)\rangle.&  \cr}$$
Now recall that $\langle q^h, z\rangle
= q^{\langle h, \eps_1 -\eps_n\rangle} =: q^\lambda$ (cf.\ 
Proposition \thtag{imageH}) and
$\langle q^{\eps_1}, z\rangle = q$.
A comparison of \eqtag{hD} and \eqtag{hC} then yields \eqtag{CD}.
\endproof
The proof of Theorem \thtag{qhC} is based on the commutation
relations \eqtag{relL} between the elements $L^\pm_{ij}\in \Uq$.
From \eqtag{deltaL} 
one derives that $L^\eps S(L^\eps) =
S(L^\eps)L^\eps = I$ ($\eps = \pm$, $I$ identity matrix).
Using this identity one can rewrite
\eqtag{relL} as
$$S(L_2^{\eps_2}) R^+ L_1^{\eps_1} = 
L_1^{\eps_1} R^+ S(L_2^{\eps_2}) \eqno\eq{SRL}$$
where $(\eps_1, \eps_2)$ is any of the three pairs
$(+,+)$, $(-,-)$, $(+,-)$. Working this out explicitly
we get:
\beginth{Lemma} lemmaA
{\thstil Suppose $1\leq i,j,k,l\leq n$
and $(\eps_1, \eps_2)\in \{(+,+), (-,-), (+,-)\}$. Then one has
the following identity:
$$\eqalignno{q^{\delta_{il}}S(L_{jl}^{\eps_2}) L_{ik}^{\eps_1} &+
(q-q^{-1}) \delta_{il}\sum_{m>i} S(L_{jm}^{\eps_2})L_{mk}^{\eps_1}
& \eq{Aa} \cr
&= q^{\delta_{kj}}L_{ik}^{\eps_1}S(L_{jl}^{\eps_2}) + (q-q^{-1})
\delta_{jk}\sum_{m<j}L_{im}^{\eps_1} S(L_{ml}^{\eps_2}).& \cr}$$
\par}
\endth\noindent
Let us introduce the following notation:
$$\eqalignno{C_1 &:= \sum_{k=2}^{n-1} q^{2(n-1)} L^+_{1k} 
S(L^-_{k1}), \quad C_2 := \sum_{k=2}^{n-1} q^{2(n-k)} L^+_{kn} 
S(L^-_{nk}), &\eq{Cidef}\cr
C_3 &:= q^{2(n-1)} L^+_{1n} S(L^-_{n1}),\quad 
C_4 := \sum_{i=1}^n q^{2(n-i)} L^+_{ii} 
S(L^-_{ii}). & \cr}$$
\beginth{Lemma} lemmaC
{\thstil One has:
$$q^hC \sim q^hC_1 + q^hC_2 + q^hC_3 + q^hC_4.\eqno\eq{Ca}$$
\par}
\endth
\beginproof
This is immediate from the definition of $C$ (cf.\ \eqtag{Cdef})
and the form of the coideal $\gok^\sigma$ (cf.\ \eqtag{kdefcd}).
\endproof\noindent
The next step is to  reduce each of the 
elements $q^hC_i$ ($1\leq i\leq 4$) modulo
$\Uq\gok^\sigma + \gok^\tau\Uq$ to a suitable element in $\Uh$. We shall do
this in a series of lemmas in which we repeatedly apply the
identities \eqtag{Aa} and \eqtag{qhLrel}. 
It will be left to the reader to fill in most of the computational
details, but with the indications we give this should be
completely straightforward.

Let us introduce 
the following terminology. Let $X$ be a monomial in the $L^\pm_{kl}$
and suppose $L^\pm_{ij}$ ($1\leq i,j \leq n$ fixed) occurs in 
\eqtag{kdefcd} (r) ($\hbox{r} = i, ii, \ldots, vi$). 
Then it is clear that
by using \eqtag{kdefcd} (r) 
we can eliminate $L^\pm_{ij}$ in the monomial $L^\pm_{ij}X$ 
modulo $\gok^\tau\Uq$. We call this procedure {\sl reduction
on the left} by \eqtag{kdefcd} (r). 
Reduction on the right is defined in a similar way. 

Observe (cf.\ \eqtag{Aa}) 
that modulo $\Uq\gok^\sigma$ we have ($2\leq k\leq n-1$):
$$ S(L^-_{k1}) L^+_{1k} 
\sim L^+_{1k}S(L^-_{k1}) + (1-q^{-2})
\sum_{m<k}L^+_{1m} S(L^-_{m1}) - (1-q^{-2}).\eqno\eq{Db}$$
Using this we can prove:
\beginth{Lemma} lemmaE
{\thstil We have ($2\leq k\leq n-1$):
$$q^h L^+_{1k}S(L^-_{k1}) \sim 
-(1-q^{-2}) \left \{ \sum_{1<m<k} q^h  L^+_{1m}S(L^-_{m1}) 
+ {1+ q^{\sigma + \tau +1 + \lambda}\over
1-q^{2\lambda + 2}}(q^{h+2\eps_1} -q^h) \right \}.$$
\par}
\endth
\beginproof
We first reduce on the left by \eqtag{kdefcd} (ii), apply
\eqtag{Aa}, and then reduce on the right by \eqtag{kdefcd} (ii). 
The result is:
$$q^h L^+_{1k}S(L^-_{k1}) \sim 
q^{-\sigma+\tau -1-\lambda} q^h S(L^-_{k1})L^+_{1k} 
+ (1-q^{-2}) q^{\tau-\lambda} 
\sum_{m<k} q^h L^-_{nm}S(L^-_{m1}).$$
Next, one eliminates $L^-_{nm}$ ($1<m<k$) and
$L^-_{n1}$ by reduction on the left by \eqtag{kdefcd}
(ii) and (vi) respectively. Then apply \eqtag{Db} to get
rid of the term $S(L^-_{k1})L^+_{1k}$:
$$\eqalign{q^h &L^+_{1k}S(L^-_{k1}) \sim \cr
&\sim q^{-\sigma+\tau -1-\lambda} q^h L^+_{1k}S(L^-_{k1})
+ (1-q^{-2}) (q^{-\sigma+\tau -1-\lambda} -1) \sum_{1<m<k}
q^h L^+_{1m}S(L^-_{m1})\cr
&+ (1-q^{-2}) (q^{-\sigma + \tau -1-\lambda} -(1-q^{2\tau}))
(q^{h+2\eps_1} -q^h) + q^{\tau +\lambda} 
(1-q^{-2}) q^h L^+_{1n} L^+_{11}.\cr}$$
Reducing $q^h L^+_{1n} L^+_{11}$
first on the left, then on the right, by \eqtag{kdefcd} (vi), 
one gets:
$$(1-q^{-2\lambda - 2})q^h L^+_{1n} L^+_{11} \sim
(q^{-\lambda} (q^{-\tau} - q^\tau) 
-q^{-2\lambda - 1} (q^{-\sigma} - q^\sigma)) (q^{h+2\eps_1}
-q^h).\eqno\eq{Ec}$$
Substituting the last formula into the preceding one,
one arrives at the desired result.
\endproof
\beginth{Lemma} lemmaF
{\thstil We have ($2\leq k\leq n-1$):
$$q^h L^+_{1k}S(L^-_{k1}) \sim
-q^{-2(k-2)}{(1-q^{-2})(1+ q^{\sigma + \tau +1 + \lambda}) 
\over 1-q^{2\lambda + 2}}(q^{h+2\eps_1} -q^h).\eqno\eq{Fa}$$
\par}
\endth
\beginproof
This follows from Lemma \thtag{lemmaE} by induction on $k$.
\endproof
In other words, we can now write $q^hC_1$ as an element of $\Uh$
modulo $\Uq\gok^\sigma + \gok^\tau\Uq$, more precisely as a linear
combination of $q^{h +2\eps_1}$ and $q^h$. We now deal with 
$q^h C_2$.
Observe (cf.\ \eqtag{Aa}) that modulo $\gok^\tau\Uq$ we have:
$$S(L^-_{nk}) L^+_{kn} + (1-q^{-2}) \sum_{m>k} S(L^-_{nm}) L^+_{mn}
\sim L^+_{kn} S(L^-_{nk}) + (1-q^{-2}).\eqno\eq{Fb}$$
\beginth{Lemma} lemmaG
{\thstil We have ($2\leq k\leq n-1$):
$$q^h L^+_{kn} S(L^-_{nk}) \sim
{(q-q^{-1})q^{-2\lambda + 1} (1 + q^{-\sigma -\tau -1 + \lambda})
\over 1-q^{-2\lambda + 2}} (q^h -q^{h-2\eps_1}) \eqno \eq{Ga}$$
\par}
\endth
\beginproof
First reduce on the left by
\eqtag{kdefcd} (iii), and then apply \eqtag{Aa}. In the resulting
expression, $L^-_{k1}$ can be eliminated by reduction on the right using
\eqtag{kdefcd} (iii). One then obtains:
$$q^h L^+_{kn} S(L^-_{nk}) \sim 
q^{\sigma-\tau+1-\lambda}q^h S(L^-_{nk})L^+_{kn}
-q^{-\tau-\lambda} (q-q^{-1}) \sum_{m>k} q^h S(L^-_{nm})L^-_{m1}.$$
Now apply \eqtag{Fb} to get rid of $S(L^-_{nk})L^+_{kn}$, 
eliminate $L^-_{m1}$ ($k<m<n$) by reducing on the right
by \eqtag{kdefcd} (iii), and replace $L^\pm_{nn}$
by $L^\mp_{11}$ (cf.\ \eqtag{kdefcd} (i)):
$$\eqalign{q^h L^+_{kn}  S(L^-_{nk})& \sim 
q^{\sigma-\tau+1-\lambda}q^h L^+_{kn} S(L^-_{nk}) \cr
& + q^{\sigma-\tau-\lambda}(q-q^{-1})(q^h -q^{h-2\eps_1}) 
- q^{-\tau-\lambda} (q-q^{-1}) q^h L^-_{11} L^-_{n1}.\cr}$$
Reducing $q^h L^-_{11} L^-_{n1}$ first on the right, then on the
left by \eqtag{kdefcd} (vi), one gets:
$$(1-q^{-2\lambda + 2}) q^h L^-_{11} L^-_{n1} \sim 
(q^{-\lambda + 1}(q^{-\tau} - q^\tau) + 
(q^\sigma - q^{-\sigma})) (q^h - q^{h-2\eps_1}).
\eqno\eq{Gb}$$
Substitution of the last formula into the preceding one yields
the lemma.
\endproof
Lemma \thtag{lemmaG} allows us to reduce $q^hC_2$ to a
linear combination of $q^h$ and $q^{h-2\eps_1}$
modulo $\Uq\gok^\sigma + \gok^\tau\Uq$.
We proceed with some preliminary results that will be used
to reduce $q^hC_3$.
\beginth{Lemma} lemmaH
{\thstil We have:
$$\sum_{m>1} S(L^-_{nm}) L^-_{m1} = \sum_{m>1}
q^{-2(m-1)+1}  L^-_{m1} S(L^-_{nm}). \eqno\eq{Ha}$$
\par}
\endth
\beginproof
We prove the following more general statement ($1\leq p\leq n$):
$$\sum_{m>1} S(L^-_{nm}) L^-_{m1} = 
\sum_{m=2}^p q^{-2(m-1)+1}  L^-_{m1} S(L^-_{nm})
+ q^{-2(p-1)} \sum_{m>p} S(L^-_{nm}) L^-_{m1}.\eqno\eq{Hb}$$
The lemma follows from this by taking $p=n$.
We prove \eqtag{Hb} by induction on $p$. The assertion is trivial for
$p=1$. Suppose the assertion is true for a fixed $p$ ($1\leq p < n$).
We use the following identity ($1\leq p < n$) 
which follows in a straightforward way
from \eqtag{Aa}:
$$\sum_{m>p} S(L^-_{nm}) L^-_{m1} = q^{-1} L^-_{p+1,1} S(L^-_{n,p+1})
+ q^{-2} \sum_{m>p+1} S(L^-_{nm}) L^-_{m1}.\eqno\eq{Hc}$$
One now proves \eqtag{Hb} with $p$ replaced by $p+1$ by first using
the induction hypothesis and then applying \eqtag{Hc}.
\endproof
\beginth{Lemma} lemmaI
{\thstil We have the following two equivalences:
$$\eqalignno{q^h L^+_{11} S(L^-_{n1}) \sim &
q^{-\tau +\lambda}\sum_{m=2}^{n-1}  q^h L^+_{1m} S(L^-_{m1}) & \eq{Ia}\cr
& \qquad \qquad
-q^{2\lambda} q^h L^+_{1n} L^+_{11}  + q^\lambda (q^{-\tau} - q^\tau)
(q^{h+2\eps_1} - q^h), & \cr}$$
$$q^h  L^-_{11} S(L^-_{n1}) \sim
q^{\tau +\lambda}\sum_{m=2}^{n-1} q^{-2(m-1)} q^h L^+_{mn} S(L^-_{nm})
-q^{-2(n-1) +1} q^h L^-_{11} L^-_{n1} \eqno\eq{Ib}$$
\par}
\endth
\beginproof
We start with the first equivalence. Reducing $q^h L^+_{11} S(L^-_{n1})$
on the left by \eqtag{kdefcd} (i) and using that 
$\sum_{m=1}^n L^-_{nm} S(L^-_{m1}) = 0$, we get:
$$q^h L^+_{11} S(L^-_{n1}) \sim
-\sum_{m=2}^{n-1} q^h L^-_{nm} S(L^-_{m1}) - q^h L^-_{n1} S(L^-_{11}).$$
Next, we eliminate $L^-_{nm}$ ($1<m<n$) and $L^-_{n1}$ by reduction
on the left using \eqtag{kdefcd} (ii) and (vi) respectively, and
simplify to obtain \eqtag{Ia}. To prove the second equivalence, observe that 
by the identity $\sum_{m=1}^n S(L^-_{nm}) L^-_{m1} = 0$ one has
$$q^h  L^-_{11} S(L^-_{n1}) \sim 
-q^{-1} \sum_{m>1} q^h S(L^-_{nm}) L^-_{m1}.$$
Apply Lemma \thtag{lemmaH} to the right-hand side,
 reduce on the left by \eqtag{kdefcd} (iii),
and use reduction on the right by \eqtag{kdefcd} (i) to eliminate
$S(L^-_{nn}) = L^+_{nn}$. One then obtains \eqtag{Ib}.
\beginth{Lemma} lemmaJ
{\thstil We have the following equivalence:
$$\eqalign{(1&-q^{-2\lambda}) q^h L^+_{1n} S(L^-_{n1}) \sim \cr
& (1-q^{-2}) \left \{ 
q^{-2\lambda} \sum_{m=2}^{n-1} q^h  L^+_{1m} S(L^-_{m1})
- q^{-\sigma + \tau -1-\lambda} 
\sum_{m=2}^{n-1} q^{-2(m-2)} q^h  L^+_{mn} S(L^-_{nm}) \right \} \cr
& + (1-q^{-2}) q^{-\sigma-2\lambda} (q^{-2(n-2)} - 1) 
q^h L^-_{11} L^-_{n1} \cr
& + q^{-\lambda} (q^{-\tau} -q^\tau) 
q^h (L^+_{11}- L^-_{11}) S(L^-_{n1}) 
-q^{-2\lambda} (q^{-\sigma} -q^\sigma) q^h S(L^-_{n1}) 
(L^+_{11}- L^-_{11}) \cr
& + (1-q^{-2}) q^{-2\lambda} (q^{h+2\eps_1} - q^{h-2\eps_1}).\cr}$$
\par}
\endth
\beginproof
To prove this lemma, we start by reducing  $q^h L^+_{1n} S(L^-_{n1})$
on the left by \eqtag{kdefcd} (vi). This leads to:
$$q^h  L^+_{1n} S(L^-_{n1}) \sim q^{-2\lambda} q^h L^-_{n1} S(L^-_{n1})
+ q^{-\lambda}(q^{-\tau} - q^\tau) q^h (L^+_{11}- L^-_{11}) S(L^-_{n1}).
\eqno\eq{Ja}$$
Now observe that by \eqtag{Aa} one has $L^-_{n1} S(L^-_{n1}) = 
S(L^-_{n1})L^-_{n1}$ and
$$\eqalignno{q S(L^-_{n1})L^+_{1n} + & (q-q^{-1}) \sum_{m>1} 
S(L^-_{nm})L^+_{mn} & \eq{Jb}\cr
& = q L^+_{1n} S(L^-_{n1}) + (q-q^{-1}) \sum_{m<n} L^+_{1m} S(L^-_{m1}).
& \cr}$$
Reducing $q^h S(L^-_{n1})  L^-_{n1}$ on the right by \eqtag{kdefcd}
(vi) and eliminating $S(L^-_{n1})L^+_{1n}$ by \eqtag{Jb}, one gets:
$$\eqalignno{q^h L^-_{n1} &S(L^-_{n1}) \sim q^h L^+_{1n} S(L^-_{n1}) &
\eq{Jc}\cr
&+ (1-q^{-2}) \sum_{m<n} q^h L^+_{1m} S(L^-_{m1}) - (1-q^{-2}) 
\sum_{m=2}^{n-1} q^h S(L^-_{nm})L^+_{mn} & \cr
& - (q^{-\sigma} -q^\sigma) q^h S(L^-_{n1})
(L^+_{11}- L^-_{11})- (1-q^{-2}) q^{h-2\eps_1}.& \cr}$$
Consider the term $\sum_{m=2}^{n-1} q^h S(L^-_{nm})L^+_{mn}$ 
in \eqtag{Jc}. We reduce
it on the right by \eqtag{kdefcd} (iii),
 apply Lemma \thtag{lemmaH}, and then
eliminate $L^-_{m1}$ ($1<m<n$) and $L^-_{n1}$ by reducing on the left by
\eqtag{kdefcd} (iii) and (vi) respectively. One ends up with:
$$\eqalign{\sum_{m=2}^{n-1} & q^h S(L^-_{nm})L^+_{mn} \sim \cr
& q^{-\sigma + \tau -1+ \lambda}\sum_{m=2}^{n-1} q^{-2(m-2)} q^h L^+_{mn}
S(L^-_{nm})  + q^{-\sigma} (1 -q^{2(n-2)}) q^h L^-_{11}L^-_{n1}.\cr}$$
Now substitute the last formula into \eqtag{Jc}, and the resulting
equation into \eqtag{Ja}. One then obtains the lemma.
\endproof
\beginth{Lemma} lemmaK
{\thstil We have:
$$q^h C_4 \sim q^{2(n-1)} q^{h+2\eps_1} + q^{h-2\eps_1}
+ {q^2-q^{2(n-1)}\over 1-q^2} q^h.\eqno\eq{Ka}$$
\par}
\endth
\beginproof
This follows by reduction using \eqtag{kdefcd} (i) and (v).
\endproof
We conclude from Lemma \thtag{lemmaC} and the subsequent lemmas
that $q^hC$ can be written modulo $\Uq \gok^\sigma + \gok^\tau \Uq$
as a linear combination of $q^{h+2\eps_1}$, $q^h$, and $q^{h-2\eps_1}$.
To conclude the proof of Theorem \thtag{qhC} we have to compute explicit
expressions for the coefficients in this linear combination.

We start with the coefficient of $q^{h+2\eps_1}$, which we denote 
by $(q^h C)^+$. Collecting terms and using
the identity 
$$\sum_{k=2}^{n-1} q^{-2(k-2)} = {1-q^{-2(n-2)}\over
1-q^{-2}},$$ we get the following expression for
$(q^hC)^+$:
$$\eqalignno{& -{(q^{2(n-1)} -q^2)
(1+q^{\sigma + \tau + 1 + \lambda})\over 
1-q^{2\lambda + 2}} - {(1-q^{-2}) (q^{2(n-1)} -q^2) q^{-2\lambda}
(1+q^{\sigma + \tau + 1 + \lambda}) \over 
(1-q^{-2\lambda})(1-q^{2\lambda + 2})}&  \cr
& - {(q^{2(n-1)} -q^2) (1+q^{\sigma + \tau + 1 + \lambda})
q^{-\tau + \lambda} M\over
(1-q^{-2\lambda}) (1-q^{2\lambda + 2})} 
- {q^{2(n-1)} q^{2\lambda} M^2 \over
(1-q^{-2\lambda}) (1-q^{-2\lambda-2})} & \cr 
& + {q^{2(n-1)} q^\lambda (q^{-\tau} -q^\tau) M\over
1-q^{-2\lambda}}  + {q^{2(n-1)} (1-q^{-2}) q^{-2\lambda}\over
1-q^{-2\lambda}} + q^{2(n-1)}, & \eq{qhCplus}\cr}$$
where
$M := q^{-\lambda}(q^{-\tau} - q^\tau) - q^{-2\lambda -1}
(q^{-\sigma} -q^\sigma)$.
Clearly, the expression in
\eqtag{qhCplus} can be uniquely written as a quotient of
a polynomial $N(q^\lambda)$ in the variable $q^\lambda$
of degree four and the polynomial 
$(1-q^{2\lambda})(1-q^{2\lambda +2})$.
A completely elementary but rather tedious calculation shows
that the coefficients of 
$N(q^\lambda) =: \sum_{i=0}^4 n_i q^{i\lambda}$
are given by:
$$\eqalign{&n_0 = 1, \quad n_1 = 
q^{\sigma + \tau + 1} + q^{-\sigma - \tau + 1} -
q^{\sigma - \tau + 1} - q^{-\sigma + \tau + 2(n-2) + 1},\cr
& n_2 = q^{2(n-2) + 2} - q^{-2\sigma + 2(n-2) + 2}
- q^{-2\tau + 2} - q^{2\tau + 2(n-2) + 2} - q^{2\sigma + 2}
+ q^2,\cr
& n_3 = q^{\sigma + \tau + 2(n-2) + 3} + q^{-\sigma - \tau + 2(n-2) + 3} 
-q^{\sigma - \tau + 3} - q^{-\sigma + \tau  + 2(n-2) + 3},\cr
& n_4 = q^{2n}.\cr}$$
From this it easily follows that $(q^h C)^+$ is equal to
$$ 
{(1+q^{\sigma + \tau + 1 + \lambda}) (1+q^{-\sigma - \tau + 1 + \lambda})
(1-q^{\sigma - \tau + 1 + \lambda})(1-q^{-\sigma + \tau + 2(n-2) +
1 + \lambda})\over (1-q^{2\lambda}) (1-q^{2\lambda + 2})},$$
in other words $(q^h C)^+ = A(q^\lambda; q^2)$ with $a,b,c,d$ as given
in Theorem \thtag{qhC}. In a similar way, one can prove that the coefficients
of $q^{h-2\eps_1}$ and $q^h$ in the reduction of $q^h C$ are as given
in Theorem \thtag{qhC}.  We shall leave the details to the reader.
However, there is a shortcut argument to prove Theorem \thtag{radialpart}
from here.
Reasoning as we did to deduce
Theorem \thtag{radialpart} from Theorem \thtag{qhC}, one  sees that
the results proved so far imply that the radial part $D$ is
a linear combination with rational coefficients in $z$
of $T_{q^2, z}$, $T_{q^{-2}, z}$ and $\id$. 
The coefficient of $T_{q^2, z}$ is equal
to $A(z; q^2)$ with parameters as given in Theorem
\thtag{radialpart}. 
Let now $w\colon \CC[z^{\pm 1}] \to \CC[z^{\pm 1}]$ be the unique algebra
automorphism sending $z$ to $z^{-1}$. The radial part $D$ commutes
with $w$, since its eigenfunctions are $w$-invariant by Proposition
\thtag{imageH}. On the
other hand, one has the relation $w\circ T_{q,z} =
T_{q^{-1}, z} \circ w$. This implies that the coefficient
of $T_{q^{-2}, z}$ in the expression for $D$ is equal to
$A(z^{-1}; q^2)$. Hence, 
the radial part $D$ and the operator $D_{AW}$ defined in \eqtag{radD}
differ at most a multiple of the identity. To show that they are actually
equal, it clearly suffices to exhibit a non-zero $P\in \CC[z+z^{-1}]$
such that $D\cdot P = D_{AW} \cdot P$. Let us take $P=1$. Then, on the
one hand, we have:
$$ D \cdot 1 = C \cdot 1 = \varepsilon(C) 1 = \sum_{i=1}^n q^{2(n-i)}\, 1
= {1-q^{2n}\over 1-q^2}\,1,$$
whereas, on the other hand, 
$D_{AW} \cdot 1 = {1-q^{2n}\over 1-q^2}\,1$ by \eqtag{radD}.
This completes the proof of Theorem \thtag{radialpart}.
\endproof

\enddocument